\documentclass[11pt]{article}

\usepackage[utf8]{inputenc}

\usepackage{amssymb,xcolor,xspace}
\usepackage{amsmath,graphicx,comment,bm}
\usepackage{gastex}

\usepackage{hyperref}
\usepackage{array}

\usepackage{geometry} 

\usepackage{fancyhdr} 
\usepackage{sectsty}
\allsectionsfont{\sffamily\mdseries\upshape} 

\usepackage[nottoc,notlof,notlot]{tocbibind} 
\usepackage[titles]{tocloft} 

\def\cqfd{\skip10=\parfillskip\parfillskip=0pt
\enspace\hfill\symbolecqfd\par\parfillskip=\skip10\par\medskip}
\def\symbolecqfd{\rlap{$\sqcap$}$\sqcup$}

\newtheorem{theorem}{Theorem}[section]
\newtheorem{proposition}[theorem]{Proposition}
\newtheorem{lemma}[theorem]{Lemma}
\newtheorem{corollary}[theorem]{Corollary}

\newtheorem{pro-fact}[theorem]{Fact}

\newtheorem{pro-example}[theorem]{Example}
\newenvironment{example}{\begin{pro-example}\rm}{\cqfd\end{pro-example}}

\newtheorem{pro-remark}[theorem]{Remark}
\newenvironment{remark}{\begin{pro-remark}\rm}{\cqfd\end{pro-remark}}

\newenvironment{preuve}{\rm \trivlist \item[\hskip \labelsep{\bf
Proof.}]}{\cqfd\endtrivlist}

\def\cqfd{\skip10=\parfillskip\parfillskip=0pt
\enspace\hfill\symbolecqfd\par\parfillskip=\skip10\par\medskip}
\def\symbolecqfd{\rlap{$\sqcap$}$\sqcup$}
\def\proof{\begin{preuve}}
\def\eop{\end{preuve}}

\def\proofof#1{\rm \trivlist \item[\hskip \labelsep{\bf
Proof of~#1.}]}
\def\eopo{\cqfd\endtrivlist}

\let\phi\varphi

\def\fffi{\textsf{fr-fi}}

\DeclareMathOperator{\relab}{norm}

\def\inv{^{-1}}
\let\epsilon\varepsilon

\def \calC {\mathcal{C}}
\def\ocalC{\overline{\mathcal{C}}}
\def \calI {\mathcal{I}}

\def \O {\mathcal{O}}

\def\calS  {\mathcal{S}}

\def\E{\mathbb{E}}

\def\MM{\mathbb{M}}

\def \P {\mathbb{P}}
\def\PSL{\textsf{PSL}_2(\Z)}

\def\Z{\mathbb{Z}}

\def\sk{\mathfrak{D}}
\def\cyc{\mathfrak{G}}
\def\rooted{\mathfrak{R}}
\def\minseq{\mathfrak{M}}
\def\probap{\mathfrak{p}}
\def\probaq{\mathfrak{q}}

\def\core{\textsf{silh}}
\def\qcore{\textsf{q-silh}}

\def\path{\textsf{path}}
\def\pathlabel{\textsf{path-label}}

\def\overlap{\mathcal{P}^{\circ\!\!\circ}}
\def\nooverlap{\mathcal{P}^{\circ\circ}}

\def \calG {\mathcal{G}}

\def\ella{\ell_2}
\def\ellb{\ell_3}
\def\ka{k_2}
\def\kb{k_3}
\def\Ella{L_2}
\def\Ellb{L_3}

\def\Kb{K_3}
\def\exc{\textsf{exc}}
\def\unroot{\textsf{unroot}}

\usepackage{todonotes}

\title{Silhouettes and generic properties of subgroups of the modular group}

\author{
    Fr\'ed\'erique Bassino, \small{\url{bassino@lipn.fr}}\\
    \small{Universit\'e Sorbonne Paris Nord, LIPN, CNRS UMR 7030, F-93430 Villetaneuse, France}%
    \and
    Cyril Nicaud, \small{\url{cyril.nicaud@u-pem.fr}}\\
    \small{LIGM, Univ Gustave Eiffel, CNRS, ESIEE Paris, F-77454, Marne-la-Vallée, France}%
    \and
    Pascal Weil, \small{\url{pascal.weil@labri.fr}}\\
    \small{Univ. Bordeaux, LaBRI, CNRS UMR 5800, F-33400 Talence, France}\thanks{%
    LaBRI, Univ. Bordeaux, 351 cours de la Lib\'eration, 33400 Talence, France.}\\
    \small{CNRS, ReLaX, UMI 2000, Siruseri, India}
    }

\begin{document}

\maketitle

\begin{abstract}
We show how to count and randomly generate finitely generated subgroups of the modular group $\PSL$ of a given isomorphism type. We also prove that almost malnormality and non-parabolicity are negligible properties for these subgroups.

The combinatorial methods developed to achieve these results bring to light a natural map, which associates with any finitely generated subgroup of $\PSL$ a graph which we call its silhouette, and which can be interpreted as a conjugacy class of free finite index subgroups of $\PSL$. 
\end{abstract}


\section{Introduction}

The modular group $\PSL$ has played a central role in algebra, number theory and geometry since the late 19th century, notably the study of its subgroups. 

A number of papers have concentrated on its finite index subgroups, with exact enumeration results for the index $n$ subgroups and results on the asymptotic behavior of that number as $n$ tends to infinity (Dey, Stothers, Muller, Schlage-Puchta and others, see \textit{e.g.} \cite{1965:Dey,1977:Stothers-Edinburgh,2004:MullerSchlage-Puchta}). Here, we deal instead with \emph{all} finitely generated subgroups of $\PSL$.

There are two main sets of results in this paper:
\begin{itemize}
\item[(R1)] enumeration and random generation results for finitely generated subgroups of $\PSL$ of a given size and isomorphism type;

\item[(R2)] proofs that almost malnormality is negligible and parabolicity is generic for finitely generated subgroups of $\PSL$.
\end{itemize}
Let us first make the notions underlying these two sets of results more explicit.

Recall that $\PSL$ is the free product of a copy of $\Z_2$ and a copy of $\Z_3$. By Kurosh's theorem (see, \textit{e.g.}, \cite{1956:Kurosh,1977:LyndonSchupp,1995:Rotman}), any finitely generated subgroup $H$ of $\PSL$ is isomorphic to a free product of $\ella$ copies of $\Z_2$, $\ellb$ copies of $\Z_3$ and $r$ copies of $\Z$: the \emph{isomorphism type} of $H$ mentioned in (R1) is the triple $(\ella, \ellb, r)$. As for results (R2), we remind the reader that $H$ is \emph{almost malnormal} if $H \cap H^x$ is finite for every $x\not\in H$; it is \emph{parabolic} if it contains a parabolic element, that is, a conjugate of a non-trivial power of $ab$, where $a$ and $b$ are the generators of the order 2 and order 3 free factors, respectively.

In addition, both sets of results refer implicitly to a distribution of probabilities on the set of finitely generated subgroups of $\PSL$, that we now explain. Each finitely generated subgroup $H$ of $\PSL$ can be represented uniquely by a finite labeled graph $\Gamma(H)$, called its \emph{Stallings graph}.
Stallings graphs, and their effective construction, were first introduced by Stallings \cite{1983:Stallings} to represent finitely generated subgroups of free groups. The idea of using finite graphs to represent subgroups of infinite, non-free groups first appeared in work of Arzhantseva and Ol'shanskii \cite{1996:ArzhantsevaOlshanskii,1998:Arzhantseva}, Gitik \cite{1996:Gitik} and Kapovich \cite{1996:Kapovich}. Markus-Epstein \cite{2007:Markus-Epstein} gave an explicit construction associating a graph with each subgroup of an amalgamated product of two finite groups, which is very close to the one used here. Here we follow the definition and construction of Kharlampovich et al. \cite{2017:KharlampovichMiasnikovWeil}. In a nutshell, the Stallings graph of a subgroup $H$ of $\PSL$ is the fragment of the Schreier (or coset) graph of $H$, spanned by the loops at vertex $H$ labeled by a geodesic representative of an element of $H$, see Section~\ref{sec: rappels}.

We take the number of vertices of $\Gamma(H)$ to be the \emph{size} of the subgroup $H$. In particular, there are only finitely many subgroups of a given size and we assume the uniform distribution on this finite set. A property of subgroups is \emph{negligible} (resp. \emph{generic}) if the proportion of size $n$ subgroups with the property tends to 0 (resp. 1) when $n$ tends to infinity.

Note that this randomness model strongly differs from those considered by Gilman \emph{et al.} in \cite{2010:GilmanMiasnikovOsin} for subgroups of hyperbolic groups and by Maher and Sisto \cite{2019:MaherSisto} for subgroups of acylindrically hyperbolic groups. Roughly speaking, these models rely on fixing an integer $k$, randomly choosing $k$ elements of $G$ using an $n$-step Markovian mechanism, considering the subgroup generated by these $k$ elements, and letting $n$ tend to infinity. It is interesting to note that in this \emph{few generators} model,  almost malnormality is generic \cite[Theorem 1.1(2)]{2019:MaherSisto}, in contrast with our model. A similar situation is already known to arise for subgroups of free groups: in the $k$-generated model, malnormality is generic (Jitsukawa \cite{2002:Jitsukawa}), whereas in the so-called graph-based model, it is negligible (Bassino \emph{et al.} \cite{2013:BassinoMartinoNicaud}). Jitsukawa's result was extended also to the Gromov density model, where the number $k$ of generators is allowed to vary as an exponential function of $n$ (Bassino \emph{et al.} \cite{2016:BassinoNicaudWeilCM}). One can argue that the model we consider in this paper is particularly natural since Stallings graphs are in bijection with subgroups.

We now return to the contributions of this paper. In \cite{2020:BassinoNicaudWeil}, the authors counted the finitely generated subgroups of $\PSL$ by size and they showed how to generate uniformly at random a subgroup of a given size. They also computed the average value of the isomorphism type of a random subgroup as a function of its size and proved a large deviations theorem for this isomorphism type. It follows that randomly generating a size $n$ subgroup of $\PSL$ will, with high probability, yield a subgroup whose isomorphism type is close to the average value. In particular, this algorithmic result does not help generate uniformly at random subgroups of a given size and isomorphism type.

In this paper, we use a completely different enumeration method for finitely generated subgroups of $\PSL$, to get a polynomial time random generation algorithm for subgroups of $\PSL$ of a given size and isomorphism type (Results (R1)). This is done by a combination of graph decomposition techniques and analytic combinatorics methods \cite{2009:FlajoletSedgewick}, involving graph operations on Stallings graphs and exponential generating series.

A construction which occurs naturally in the proofs of these results is what we call the \emph{silhouetting} of the Stallings graph of a finitely generated subgroup of $\PSL$. It consists in a sequence of ``simplifications'' of the graph, leading (except in extremal cases) to a uniform degree loop-free graph, which represents a conjugacy class of finite index, free subgroups of $\PSL$. A remarkable, and somewhat surprising property is that the silhouetting construction preserves uniformity. More precisely, among the size $n$ subgroups whose silhouette has size $s$, every size $s$ silhouette graph is equally likely.

This property allows us to lift asymptotic properties of silhouette graphs, which are more easily understood, to all Stallings graphs of finitely generated subgroups of $\PSL$. This leads to Results (R2), namely the negligibility of almost malnormality and the genericity of parabolicity. The proof of these results is of a combinatorial nature, and rather intricate.

We would like to point out in particular an intermediate result which may be of independent interest. We show that, with high probability, in a finite group of permutations generated by a pair of fixpoint-free permutations $(\sigma_2,\sigma_3)$, of order 2 and 3 respectively, the composition $\sigma_2\sigma_3$ admits orbits of a certain, relatively small size (Proposition~\ref{pro:small cycles}). Obtaining such results on the composition of two randomly chosen mappings is notoriously difficult. It is, for instance, a bottleneck in the study of the properties of random deterministic automata~\cite{2014:Nicaud}. Most known results rely on a fine grain independent analysis of the mappings, but we know very little on their composition. Character theory has been used to tackle this kind of difficulties in the study of combinatorial maps~\cite{1990:JacksonVisentin,2006:Gamburd}: this approach yields enumeration results on triplets $(\sigma_1,\sigma_2,\tau)$ such that $\sigma_1\circ\sigma_2=\tau$, but for a fixed cyclic type of $\tau$ only, and it seems very difficult to exploit such results for our purposes (see also \cite{2019:BudzinskiCurienPetri}).

We note that, while there are a good number of results in the literature about the genericity or negligibility of certain properties of subgroups of free groups (\cite{2002:Jitsukawa,2013:BassinoMartinoNicaud,2008:BassinoNicaudWeil,2016:BassinoNicaudWeilCM,2016:BassinoNicaudWeil}), there are precious few such results for subgroups of other groups. We can cite in this direction, for their pioneering methods, the results of Arzhantseva and Ol'shanskii \cite{1996:ArzhantsevaOlshanskii} who show for instance that, for a very large (generic) class of $r$-generator, $k$-relator presentations ($r,k$ fixed), all $\ell$-generated subgroups ($\ell < r$) are free and quasi-convex; and the results mentioned above of Gilman \emph{et al.} \cite{2010:GilmanMiasnikovOsin} and Maher and Sisto \cite{2019:MaherSisto} on $k$-generated subgroups of a fixed hyperbolic or acylindrically hyperbolic group. To our knowledge, our results (R2) on almost malnormality  and parabolicity for subgroups of the modular group are the first that are based on the distribution of subgroups given by Stallings graphs.

The paper is organized as follows. Readers can find in Section~\ref{sec: rappels} the precise definitions of the Stallings graph of a subgroup of $\PSL$ and its combinatorial type, and results from the literature relating this combinatorial information with algebraic properties of the subgroup such as its isomorphism type, its index or its freeness.

Section~\ref{sec: moves and silhouette} introduces combinatorial operations on Stallings graphs. Iterating these operations is a confluent process, which leads to the so-called \emph{silhouette} of the given graph or subgroup. We show there that silhouetting preserves uniformity (Theorems~\ref{thm:uniformity} and~\ref{thm:uniformity rooted}), and that the size of the silhouette of a size $n$ subgroup is super-polynomially generically at least $n-3n^{\frac23}$(Propositions~\ref{sec: size silhouette} and~\ref{prop: large silhouette rooted}). Also, silhouetting preserves the free rank component of the isomorphism type of a subgroup (Proposition~\ref{prop: rank and small silhouette}).

In Section~\ref{sec: counting}, we exploit a fine description of the operations defined in Section~\ref{sec: moves and silhouette} to give exact counting formulas for the number of subgroups of $\PSL$ of a given combinatorial or isomorphism type. A key ingredient here, namely the pointing operation on exponential generating series, is borrowed from the toolbox of analytic combinatorics \cite{2009:FlajoletSedgewick}.

Section~\ref{sec: random generation} uses the operations from Section~\ref{sec: moves and silhouette} in a different way to formulate an algorithm (which includes a rejection algorithm component) to generate uniformly at random a subgroup of a given size and isomorphism type.

The last section, Section~\ref{sec: generic properties}, shows that, generically, subgroups of $\PSL$ contain parabolic elements (Corollary~\ref{cor: parabolicity}) and fail to be almost malnormal (Theorem~\ref{thm: negligibility}). Both results exploit the statistical result on the existence of cycles labeled by non-trivial powers of $ab$ (where $a$ and $b$ are the order 2 and order 3 generators of $\PSL$) in the Stallings graph of a subgroup mentioned above (Proposition~\ref{pro:small cycles}). 

To conclude this introductory section, we note that, according to the results presented here, the silhouetting operation is combinatorially and asymptotically significant in the study of finitely generated subgroups of $\PSL$. We are able to lift a statistical property of silhouette graphs to the class of all $\PSL$-reduced graphs, and it would be interesting to see what other properties can be lifted in that fashion. In addition, we think that the silhouetting operation also has a topological, or possibly a geometric interpretation, even for finite index subgroups, and we would be curious about its properties.

\section{Preliminaries}\label{sec: rappels}

We work with the following presentation of the modular group:
$$\PSL = \langle a,b \mid a^2 = b^3 = 1\rangle.$$
The elements of $\PSL$ are represented by words over the alphabet $\{a,b,a\inv,b\inv\}$. Since $a\inv = a$ in $\PSL$, we can eliminate the letter $a\inv$ from this alphabet. Each element of $\PSL$ then has a unique shortest (or \emph{normal}, or \emph{geodesic}) representative, which is a freely reduced word without factors in $\{a^2, b^2, b^{-2}\}$. That is, the normal representatives are the words of length at most 1 and the words alternating letters $a$ and letters in $\{b,b\inv\}$.

The \emph{Schreier graph} (or \emph{coset graph}) of a subgroup $H$ of $\PSL$ is the graph whose vertices are the cosets $Hg$ of $H$ ($g\in \PSL$), with an $a$-labeled edge from $Hg$ to $Hga$ and a $b$-labeled edge from $Hg$ to $Hgb$, for every $g\in G$. We think of $b$-edges as 2-way edges, reading $b$ in the forward direction and $b\inv$ in the backward direction. Since $a = a^2$, there is an $a$-edge from vertex $v$ to vertex $v'$ if and only there is one from $v'$ to $v$: as a result, we think of the $a$-edges as undirected edges, that can be traveled in either direction, each time reading $a$. 

Note that a word is in $H$ if and only if it labels a loop at vertex $v_0 = H$ in the Schreier graph of $H$. The \emph{Stallings graph} of $H$, written $(\Gamma(H),v_0)$, is the fragment of the Schreier graph of $H$ spanned by the loops at $v_0$ reading the geodesic representatives of the elements of $H$, rooted at $v_0$. In particular, a word is in $H$ if and only if its geodesic representative labels a loop in $\Gamma(H)$ at vertex $v_0$. We refer the reader to Remark~\ref{rk: historical} and to \cite{2017:KharlampovichMiasnikovWeil} for more details on these graphs. We note in particular that $H$ has a finite Stallings graph if and only if it is finitely generated, and that $\Gamma(H)$ is efficiently algorithmically computable if $H$ is given by a finite set of generators (words on the alphabet $\{a,b,b\inv\}$) \cite{2017:KharlampovichMiasnikovWeil,2020:BassinoNicaudWeil}.

\begin{example}\label{ex: Stallings graphs}
Figure~\ref{fig: stallings graphs} shows examples of Stallings graphs. 
\begin{figure}[htbp]
\centering
\begin{picture}(100,100)(0,-98)
\gasset{Nw=4,Nh=4}
\node(n0)(0.0,-0.0){$1$}

\node(n1)(6.0,-14.0){2}
\node(n2)(0.0,-28.0){3}
\node(n3)(34,-0.0){4}
\node(n4)(28,-14.0){5}
\node(n5)(34.0,-28.0){6}

\drawedge(n0,n1){$b$}

\drawedge(n1,n2){$b$}

\drawedge(n2,n0){$b$}

\drawedge[ELside=r](n3,n4){$b$}

\drawedge[ELside=r](n4,n5){$b$}

\drawedge[ELside=r](n5,n3){$b$}

\drawedge[AHnb=0](n0,n3){$a$}

\drawedge[AHnb=0](n1,n4){$a$}

\drawedge[AHnb=0](n2,n5){$a$}
\node(m1)(72.0,-4.0){$1$}
\node(m2)(78.0,-14.0){2}
\node(m3)(72.0,-24.0){3}
\node(m4)(52,-4.0){4}
\node(m5)(52,-24.0){5}
\node(m6)(98.0,-14.0){6}

\drawedge(m1,m2){$b$}

\drawedge(m2,m3){$b$}

\drawedge(m3,m1){$b$}

\drawedge[AHnb=0](m1,m4){$a$}

\drawedge[ELside=r](m4,m5){$b$}

\drawedge[AHnb=0](m5,m3){$a$}

\drawedge[AHnb=0](m2,m6){$a$}

\drawloop(m6){$b$}
\node(h11)(10.0,-40.0){$1$}
\node(h44)(34.0,-40.0){2}
\node(h55)(58.0,-40.0){3}
\node(h23)(82.0,-40.0){4}

\node(h22)(10.0,-68.0){5}
\node(h33)(2.0,-54.0){6}
\node(h12)(82.0,-68.0){7}
\node(h31)(90.0,-54.0){8}

\node(h65)(2.0,-76.0){9}
\node(h54)(16.0,-84.0){10}
\node(h70)(30.0,-76.0){11}

\node(h71)(30.0,-96.0){12}

\node(h56)(34.0,-68.0){13}
\node(h46)(58.0,-68.0){14}
\node(h66)(46.0,-76.0){15}

\node(h32)(62.0,-88.0){16}
\node(h64)(98,-74.0){17}
\node(h45)(98.0,-96.0){18}
\node(h13)(74.0,-96.0){19}
\node(h21)(50.0,-96.0){20}

\drawedge(h11,h22){$b$}
\drawedge[ELside=r](h22,h33){$b$}
\drawedge(h33,h11){$b$}

\drawedge(h23,h31){$b$}
\drawedge[ELside=r](h31,h12){$b$}
\drawedge(h12,h23){$b$}

\drawedge(h13,h21){$b$}
\drawedge(h21,h32){$b$}
\drawedge(h32,h13){$b$}

\drawedge[ELside=r](h46,h56){$b$}
\drawedge[ELside=r](h56,h66){$b$}
\drawedge[ELside=r](h66,h46){$b$}

\drawedge(h54,h65){$b$}
\drawedge(h65,h70){$b$}
\drawedge(h70,h54){$b$}

\drawedge[AHnb=0](h22,h56){$a$}

\drawedge[AHnb=0](h12,h46){$a$}

\drawedge[AHnb=0](h32,h66){$a$}

\drawedge[AHnb=0](h11,h44){$a$}
\drawedge(h44,h55){$b$}
\drawedge[AHnb=0](h55,h23){$a$}

\drawedge[AHnb=0](h31,h64){$a$}
\drawedge(h64,h45){$b$}
\drawedge[AHnb=0](h45,h13){$a$}

\drawedge[AHnb=0](h70,h71){$a$}

\drawloop[loopangle=-180](h71){$b$}
\drawloop[loopangle=-135,AHnb=0](h54){$a$}
\drawloop[loopangle=180,AHnb=0](h21){$a$}

\drawedge[AHnb=0](h65,h33){$a$}
\end{picture}
\caption{Top: the Stallings graphs of the subgroups
$H = \langle abab\inv, babab\rangle$ and $K = \langle abab, babab\inv\rangle$ of $\PSL$. Bottom: the Stallings graph of 
  $L = \langle (b\inv a)^2bab, b\inv (ab)^2ab\inv ab, (abab\inv)^2, (ba)^2b\inv a(ba)^3b\inv a,(ab)^8ab\inv a\rangle$.
In each case, the root is the vertex labeled 1.}\label{fig: stallings graphs}
\end{figure}
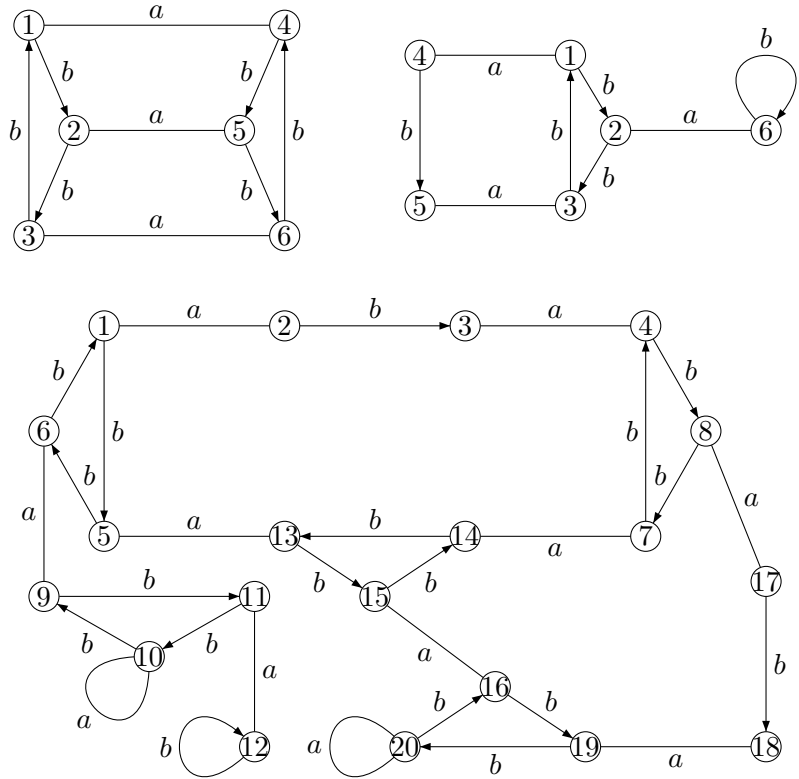
\end{example}

\begin{remark}\label{rk: historical}
The definition of Stallings graphs given above is a generalization of that introduced by Stallings in 1983 \cite{1983:Stallings} for finitely generated subgroups of free groups, and a particular instance of the definition first introduced by Gitik \cite{1996:Gitik} in 1996 under the name of geodesic core, and systematized by Kharlampovich \textit{et al.} \cite{2017:KharlampovichMiasnikovWeil} in 2017. Given a finitely presented group $G = \langle A\mid R\rangle$, a language $L$ of representatives for $G$ (a set of words over the alphabet $A \cup A\inv$) and a subgroup $H$, one considers the fragment of the Schreier graph of $H$ spanned by the $L$-representatives of the elements of $H$. It is effectively computable if $G$ is equipped with an automatic structure \cite{2017:KharlampovichMiasnikovWeil}. In the particular case where $G = \PSL$, the algorithm is quite straightforward, we refer the reader to \cite{2020:BassinoNicaudWeil} for an outline.

This notion is not the first usage of graphs (or automata-theoretic representations, see Gersten and Short \cite{1991:GerstenShort,1991:Short}) to reason about subgroups of infinite groups: in 1996, Arzhantseva and Ol'shanskii \cite{1996:ArzhantsevaOlshanskii} used graphs with great success to investigate the subgroups of an exponentially generic class of $k$-relator groups; the same year, Kapovich \cite{1996:Kapovich} used a different approach to compute the quasi-convexity constant of subgroups of automatic groups. However, the systematic approach and the unique representation provided by Stallings graphs are extremely fruitful in tackling decision and computability problems, see \cite{2017:KharlampovichMiasnikovWeil}, or enumeration problems \cite{2020:BassinoNicaudWeil}. 

The particular case of interest here, namely the free product of two finite groups, was investigated from an algorithmic point of view in 2007 by Markus-Epstein \cite{2007:Markus-Epstein}. The graphs she computes are different from (and a little larger than) the ones we consider here, because the language $L$ of representatives she uses strictly contains the language of geodesics.
\end{remark}

It is immediate from the definition of Stallings graphs that the $a$-edges (resp. $b$-edges) of $\Gamma(H)$ form a partial, injective map on the vertex set of the graph. Moreover, because $a^2 = b^3 = 1$, distinct $a$-edges are never adjacent to the same vertex: we distinguish therefore $a$-loops and so-called \emph{isolated $a$-edges}. Similarly, if we have two consecutive $b$-edges, say from $v_1$ to $v_2$ and from $v_2$ to $v_3$, then $\Gamma(H)$ also has a $b$-edge from $v_3$ to $v_1$. Thus the $b$-edges are either loops, or \emph{isolated $b$-edges}, or part of a $b$-triangle. Finally, every vertex except maybe the root vertex is adjacent to an $a$- and to a $b$-edge.

A rooted graph satisfying these conditions is called \emph{$\PSL$-reduced} and it is not difficult to see that every finite $\PSL$-reduced rooted graph is the Stallings graph of a unique finitely generated subgroup of $\PSL$. That is, the mapping $H \mapsto (\Gamma(H),v_0)$ is a bijection between finitely generated subgroups of $\PSL$ and $\PSL$-reduced rooted graphs.

The \emph{combinatorial type} of a $\PSL$-reduced graph $\Gamma$ is the tuple $(n,\ka,\kb,\ella,\ellb)$ where $n$ is the number of vertices of $\Gamma$, $\ka$ and $\kb$ are the numbers of isolated $a$- and $b$-edges, and $\ella$ and $\ellb$ are the numbers of $a$- and $b$-loops. We sometimes talk of the combinatorial type of a subgroup to mean the combinatorial type of its Stallings graph, and we refer to $n$ (the number of vertices of a $\PSL$-reduced graph) as the \emph{size} of the graph or even the size of the subgroup. See \cite{2020:BassinoNicaudWeil} for a discussion of the possible combinatorial types. 

An edge-labeled graph $\Gamma$ is said to be \emph{$\PSL$-cyclically reduced} if it is $\PSL$-reduced when rooted at every one of its vertices. We also say that a finitely generated subgroup of $\PSL$ is $\PSL$-cyclically reduced if its Stallings graph is. In that case, the combinatorial type $(n,\ka,\kb,\ella,\ellb)$ satisfies the following constraints:
$n = 2\ka + \ella = 2\kb + \ellb + 3m$ where $m$ is the number of $b$-triangles.

Let us also record the following results (see, \emph{e.g.}, \cite{2020:BassinoNicaudWeil}).

\begin{proposition}\label{prop: charact free and findex}
A subgroup $H\le \PSL$ has finite index if and only if its Stallings graph is $\PSL$-cyclically reduced and has combinatorial type of the form $(n,\ka,0,\ella,\ellb)$. It is free if and only if its combinatorial type is of the form $(n,\ka,\kb,0,0)$.

Proper free $\PSL$-cyclically reduced subgroups have even size, and proper free and finite index subgroups have size a positive multiple of 6.
\end{proposition}

By Kurosh's classical theorem on subgroups of free groups (\textit{e.g.}, \cite[Proposition III.3.6]{1977:LyndonSchupp}), a subgroup $H$ of $\PSL$ is isomorphic to a free product of $r_2$ copies of $\Z_2$, $r_3$ copies of $\Z_3$ and a free group of rank $r$, for some non-negative integers $r_2,r_3,r$. The triple $(r_2,r_3,r)$, which characterizes $H$ up to isomorphism (but not up to an automorphism of $\PSL$) is called the \emph{isomorphism type} of $H$.
We record the following connection between the combinatorial and the isomorphism types of a subgroup \cite[Proposition 2.9]{2020:BassinoNicaudWeil}.

\begin{proposition}\label{prop: combinatorial vs. isomorphism}
Let $H$ be a subgroup of $\PSL$ of size at least 2 and let $(n,\ka,\kb,\ella,\ellb)$ be the combinatorial type of $\Gamma(H)$.

If $\Gamma(H)$ is $\PSL$-cyclically reduced, the isomorphism type of $H$ is
$$(\ella,\ellb,1+\frac{n-2\kb-3\ella-4\ellb}6).$$

If $\Gamma(H)$ is not $\PSL$-cyclically reduced, the isomorphism type of $H$ is
\begin{align*}
(\ella,\ellb,\frac13+\frac{n-2\kb-3\ella-4\ellb}6) & \quad\text{if the base vertex is adjacent to an $a$-edge} \\
(\ella,\ellb,\frac23+\frac{n-2\kb-3\ella-4\ellb}6) & \quad\text{if the base vertex is adjacent to a $b$-edge.}
\end{align*}
\end{proposition}

One of our objectives in this paper is to count subgroups by isomorphism type or by combinatorial type. Since subgroups are in bijection with $\PSL$-reduced rooted graphs (their Stallings graphs), it is equivalent to count these graphs. For technical reasons, it is easier to count \emph{labeled graphs}, that is, graphs whose vertex set is equipped with a (labeling) bijection onto a set of the form $[n] = \{1,\dots,n\}$. The graphs in Figure~\ref{fig: stallings graphs} are in fact labeled graphs.

The prefered tool to count and randomly generate labeled structures (graphs or rooted graphs in this paper) is the theory of exponential generating series (EGS), to which we can apply the powerful tools of analytic combinatorics, see \cite{2009:FlajoletSedgewick}. More details are introduced where they are needed.

\begin{example}\label{ex: small n}
The $\PSL$-cyclically reduced graphs $\Gamma$ with 1 or 2 vertices are represented in Figure~\ref{fig: 2-vertex cyclically reduced}. 
There is only one with 1 vertex, and three with 2 vertices. Two of them can be labeled in two different ways while the third one admits only one labeling.
\end{example}
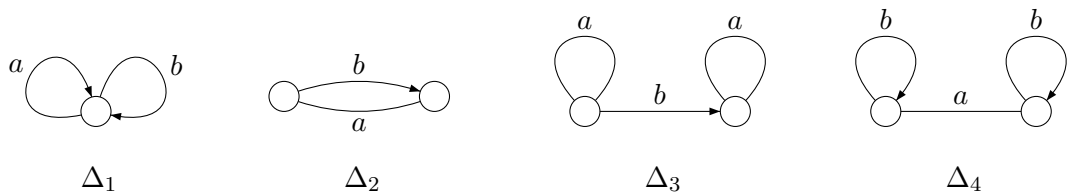
\begin{figure}[htbp]
\centering
\begin{picture}(135,19)(3,-17)
\gasset{Nw=4,Nh=4}

\node(n0)(8.0,-5.0){}
\put(6,-15){$\Delta_1$}
\drawloop[loopangle=150](n0){$a$}
\drawloop[loopangle=30](n0){$b$}

\node(n0)(33.0,-3.0){}
\node(n1)(53.0,-3.0){}
\put(41,-15){$\Delta_2$}
\drawedge[curvedepth=2](n0,n1){$b$}
\drawedge[curvedepth=2,AHnb=0](n1,n0){$a$}

\node(n0)(73.0,-5.0){}
\node(n1)(93.0,-5.0){}
\put(81,-15){$\Delta_3$}
\drawedge(n0,n1){$b$}
\drawloop[AHnb=0,loopangle=90](n0){$a$}
\drawloop[AHnb=0,loopangle=90](n1){$a$}

\node(n0)(113.0,-5.0){}
\node(n1)(133.0,-5.0){}
\put(121,-15){$\Delta_4$}
\drawedge[AHnb=0](n0,n1){$a$}
\drawloop[loopangle=90](n0){$b$}
\drawloop[loopangle=90](n1){$b$}
\end{picture}
\caption{\small All $\PSL$-cyclically reduced graphs with at most 2 vertices.}\label{fig: 2-vertex cyclically reduced}
\end{figure}

\section{Moves on $\PSL$-reduced graphs and silhouette of a subgroup}\label{sec: moves and silhouette}

We will see in Section~\ref{sec: reduction to cyclically reduced} below that counting and randomly generating subgroups of $\PSL$ reduces to counting and randomly generating labeled $\PSL$-cyclically reduced graphs. 
Before we embark on this task, we introduce a combinatorial construction on this class of graphs.

More precisely, we define in Section~\ref{sec: moves} a number of moves on a labeled $\PSL$-cyclically reduced graph, depending on its geometry. They are used in Section~\ref{sec: silhouette} to identify an interesting structure within a $\PSL$-cyclically reduced graph, namely its \emph{silhouette}. They will later be used in Section~\ref{sec: couting cyclically reduced} to count subgroups of $\PSL$ and in Section~\ref{sec: random generation} to randomly generate them.

Sections~\ref{sec: randomness preserved} and~\ref{sec: size silhouette} discuss technical properties of the silhouetting operation, which will be used in Section~\ref{sec: generic properties} to establish certain asymptotic properties of subgroups of $\PSL$.

Finally, we examine in Section~\ref{sec: silhouette all subgroups} how the silhouetting operation can be extended to all subgroups of $\PSL$, even if not $\PSL$-cyclically reduced.

\subsection{Moves on a labeled $\PSL$-cyclically reduced graph}\label{sec: moves}

Here we define so-called $\lambda_3$-, $\lambda_{2,1}$-, $\lambda_{2,2}$-, $\kappa_3$- and $\exc$-moves on quasi-labeled $\PSL$-cyclically reduced graphs (see below for the definition of quasi-labeled structures). 

Each of these moves deletes vertices from the input graph, so that the resulting graph has an injective labeling function (defined on the vertex set) which is into, but possibly not onto, a set of the form $[n]$. We then talk of a \emph{quasi-labeled} graph. It is clear that a quasi-labeled graph $\Gamma$ can be canonically relabeled (by a uniquely defined order-preserving map) into a properly labeled graph, which we write $\relab(\Gamma)$.

If $\Gamma$ is a quasi-labeled $\PSL$-cyclically reduced graph, we also introduce the notion of \emph{rank} of an isolated $a$-edge. Note that these edges are totally ordered by comparing the minimum label of one of their end vertices. We say that an isolated $a$-edge $e$ has rank $i \ge 1$ if there are exactly $i-1$ isolated $a$-edges that are less than $e$.

Let $\Gamma$ be a quasi-labeled $\PSL$-cyclically reduced graph, with combinatorial type $\bm\tau$. To lighten up the description of the moves we define on $\Gamma$, we abusively identify the vertices of $\Gamma$ with their labels. We also abusively write $\Delta_i$ ($i = 1,2,3,4$) for any quasi-labeled version of the graphs in Example~\ref{ex: small n}.

\paragraph{$\lambda_3$-moves}
Let $\ell$ be a $b$-loop in $\Gamma$, say at vertex $v$ (in fact, at the vertex labeled $v$). Then we are in exactly one of the following situations.
\begin{itemize}
\item[(1)] $v$ is adjacent to an isolated $a$-edge, connecting it with a vertex $w\ne v$.

\item[(2)] vertex $v$ carries an $a$-loop; in that case, $\Gamma = \Delta_1$.
\end{itemize}

\begin{figure}[htbp]
\centering
\begin{picture}(130,20)(-2,-22)
\gasset{Nw=4,Nh=4}
\drawoval(37,-12,40,20,12)
\node(n0)(20.0,-12.0){$w$}
\node(n1)(4.0,-12.0){$v$}
\drawedge[AHnb=0](n0,n1){$a$}
\drawloop[loopangle=180,loopdiam=6](n1){$b$}

\drawoval(110,-12,40,20,12)
\node(n2)(93.0,-12.0){$w$}
\drawloop[loopangle=180,loopdiam=6,AHnb=0](n2){$a$}
\end{picture}
\caption{\small The graphs $\Gamma$ and $\Delta$ in Case (1)}\label{fig: tb differenciation}
\end{figure}
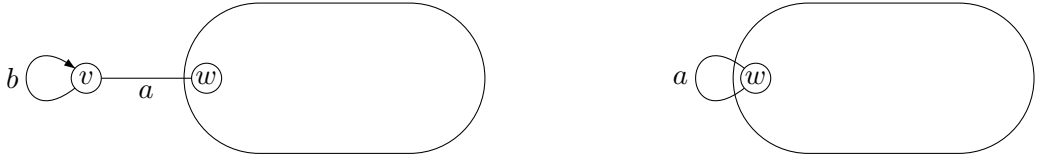

In Case (1), the \emph{$(\lambda_3,v,w)$-move} consists in deleting vertex $v$ and the adjacent edges, and adding an $a$-loop $\ell'$ at vertex $w$. The resulting quasi-labeled graph $\Delta$ (see Figure~\ref{fig: tb differenciation}) is $\PSL$-cyclically reduced graph  and has combinatorial type $\bm\tau + \bm\lambda_3$, where $\bm\lambda_3 =  (-1,-1,0,1,-1)$.

Note that the pair $(\Gamma,\ell)$ can be retrieved unequivocally from $\Delta$ and the triple $(\lambda_3,v,w)$ --- provided that the quasi-labeled graph $\Delta$ has no vertex $v$, has a vertex $w$ and has an $a$-loop at $w$, we say that $\Delta$ is \emph{valid for $(\lambda_3,v,w)$}.

Switching to labeled graphs (namely $\relab(\Gamma)$ and $\relab(\Delta)$), we note that $\lambda_3$-moves in general establish a bijection between the set of structures $(\Gamma,\ell)$ that arise from Case (1) with $\Gamma$ (properly) labeled, of size $n$ and combinatorial type $\bm\tau$, and the set of structures $(\Delta,\ell',v)$ formed by a labeled $\PSL$-cyclically reduced graph $\Delta$ with combinatorial type $\bm\tau + \bm\lambda_3$, an $a$-loop $\ell'$ in $\Delta$ and an integer $v\in [n]$.

No move is defined in Case (2).

\paragraph{$\lambda_2$-moves}
Let $\ell$ be an $a$-loop in $\Gamma$, say at vertex $v$. Then we are in exactly one of the following situations.
\begin{itemize}
\item[(1)] vertex $v$ sits on a $b$-triangle (equivalently: it does not carry a $b$-loop, and it has an incoming and an outgoing $b$-edge).

\item[(2)] vertex $v$ is adjacent to an isolated $b$-edge, linking it with a vertex $w\ne v$, which in turn is adjacent to an isolated $a$-edge linking it with a third vertex $w'\ne v,w$.

\item[(3)] vertex $v$ is adjacent to an isolated $b$-edge, linking it with a vertex $w\ne v$, which carries an $a$-loop; in that case $\Gamma = \Delta_3$.

\item[(4)] vertex $v$ carries a $b$-loop; in that case, $\Gamma = \Delta_1$.
\end{itemize}

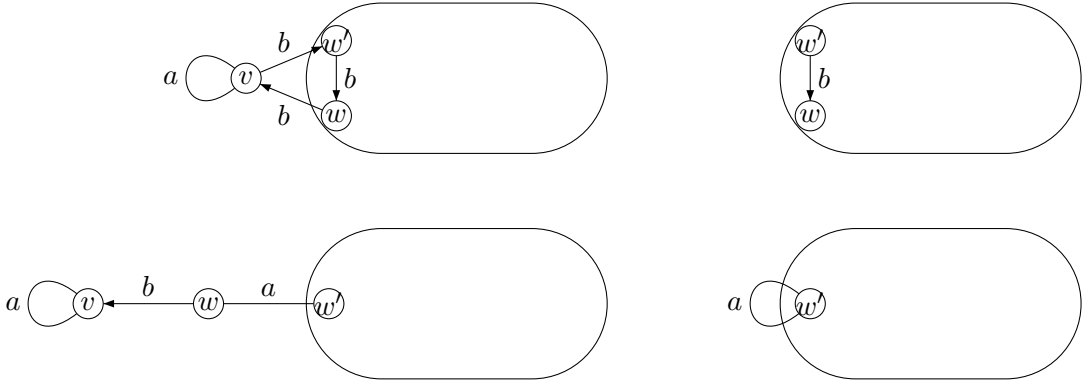
\begin{figure}[htbp]
\centering
\begin{picture}(130,50)(-2,-52)
\gasset{Nw=4,Nh=4}
\drawoval(49,-12,40,20,12)
\node(n0)(33.0,-7.0){$w'$}
\node(n00)(33.0,-17.0){$w$}
\node(n1)(21.0,-12.0){$v$}
\drawedge(n00,n1){$b$}
\drawedge(n1,n0){$b$}
\drawedge(n0,n00){$b$}
\drawloop[loopangle=180,loopdiam=6,AHnb=0](n1){$a$}

\drawoval(112,-12,40,20,12)
\node(n2)(96.0,-7.0){$w'$}
\node(n22)(96.0,-17.0){$w$}
\drawedge(n2,n22){$b$}

\drawoval(49,-42,40,20,12)
\node(n3)(32.0,-42.0){$w'$}
\node(n4)(16.0,-42.0){$w$}
\node(n5)(0.0,-42.0){$v$}
\drawedge[ELside=r](n4,n5){$b$}
\drawedge[ELside=r,AHnb=0](n3,n4){$a$}
\drawloop[loopangle=180,loopdiam=6,AHnb=0](n5){$a$}

\drawoval(112,-42,40,20,12)
\node(n6)(96.0,-42.0){$w'$}
\drawloop[loopangle=180,loopdiam=6,AHnb=0](n6){$a$}
\end{picture}
\caption{\small The labeled graphs $\Gamma$ and $\Delta$ in Cases (1) and (2)}\label{fig: ta differenciation}
\end{figure}

In Case (1), let $w$ and $w'$ be the other extremities of the $b$-edges ending and starting at $v$, respectively. Then $w\ne w'$ and $\Gamma$ has a (non-isolated) $b$-edge from $w$ to $w'$. The $(\lambda_{2,1},v,w')$-move consists in removing from $\Gamma$ vertex $v$ and the adjacent edges (the $a$-loop $\ell$ and two $b$-edges). The resulting graph $\Delta$ (see Figure~\ref{fig: ta differenciation}) is $\PSL$-cyclically reduced, it has an isolated $b$-edge from $w'$ to $w$ and combinatorial type $\bm\tau  + \bm\lambda_{2,1}$, where $\bm\lambda_{2,1} =  (-1,0,1,-1,0)$.

Here too, the pair $(\Gamma,\ell)$ can be retrieved unequivocally from $\Delta$ and the triple $(\lambda_{2,1},v,w')$ --- provided that the quasi-labeled graph $\Delta$ has no vertex $v$, has a vertex $w'$ and has an isolated $b$-edge starting at $w'$, we say that $\Delta$ is \emph{valid for $(\lambda_{2,1},v,w')$}.

In particular, $\lambda_{2,1}$-moves establish a bijection between the set of structures $(\Gamma,\ell)$ that arise from Case (1) with $\Gamma$ (properly) labeled, of size $n$ and combinatorial type $\bm\tau$, and the set of structures $(\Delta,e',v)$ formed by a labeled $\PSL$-cyclically reduced graph $\Delta$ with combinatorial type $\bm\tau + \bm\lambda_{2,1}$, an isolated $b$-edge $e'$ in $\Delta$ and an integer $v\in [n]$.

\medskip

In Case (2), the \emph{$(\lambda_{2,2}, v\to w,w')$-move} (resp. $(\lambda_{2,2}, v\leftarrow w,w')$, depending on the orientation of the $b$-edge adjacent to $v$) consists in deleting from $\Gamma$ the vertices $v$ and $w$ and the edges adjacent to them, and adding an $a$-loop $\ell'$ at $w'$. The resulting graph $\Delta$ (see Figure~\ref{fig: ta differenciation}) is $\PSL$-cyclically reduced and has combinatorial type $\bm\tau  + \bm\lambda_{2,2}$, where $\bm\lambda_{2,2} =  (-2,-1,-1,0,0)$.

Here, the pair $(\Delta,\ell')$ and the triple $(\lambda_{2,2}, v\to w,w')$ (resp. $(\lambda_{2,2}, v\leftarrow w,w')$) specify uniquely the input pair $(\Gamma,\ell)$ -- provided that $\Delta$ is \emph{valid for $(\lambda_{2,2}, v,w,w')$}, namely that $\Delta$ has no vertices $v$ or $w$ and has an $a$-loop at a vertex $w'$. Note that this validity allows reversing both a $(\lambda_{2,2}, v\to w,w')$-move and a $(\lambda_{2,2}, v\leftarrow w,w')$-move, yielding different quasi-labeled graphs $\Gamma$.

It follows that $\lambda_{2,2}$-moves establish a bijection between the set of structures $(\Gamma,\ell)$ that arise from Case (2) with $\Gamma$ (properly) labeled, of size $n$ and combinatorial type $\bm\tau$, and the set of structures $(\Delta,\ell',v,w,\epsilon)$ formed by a labeled $\PSL$-cyclically reduced graph $\Delta$ with combinatorial type $\bm\tau + \bm\lambda_{2,2}$, an $a$-loop $\ell'$ in $\Delta$, distinct integers $v,w\in [n]$ and a boolean $\epsilon$ (to account for the choice between a $(\lambda_{2,2}, v\to w,w')$- and a $(\lambda_{2,2}, v\leftarrow w,w')$-move).

No $\lambda_2$-move is defined in Cases (3) and (4).

\paragraph{$\kappa_3$-moves}
Let $e$ be an isolated $b$-edge in $\Gamma$, say from vertex $v$ to vertex $w \ne v$. Then we are in exactly one of the following situations.
\begin{itemize}
\item[(1)] Vertices $v$ and $w$ are adjacent to isolated $a$-edges, linking them respectively to vertices $v'$ and $w'$, and $v,w,v',w'$ are pairwise distinct.

\item[(2)] Exactly one of the $v$ and $w$ is adjacent to an $a$-loop.

\item[(3)] Vertices $v$ and $w$ are adjacent to the same isolated $a$-edge; in this case $\Gamma = \Delta_2$.

\item[(4)] Vertices $v$ and $w$ are both adjacent to an $a$-loop; in that case, $\Gamma = \Delta_3$.
\end{itemize}

\begin{figure}[htbp]
\centering
\begin{picture}(130,20)(-2,-22)
\gasset{Nw=4,Nh=4}
\drawoval(37,-12,40,20,12)
\node(n0)(21.0,-07.0){$v'$}
\node(n00)(21.0,-17.0){$w'$}
\node(n1)(4.0,-07.0){$v$}
\node(n11)(4.0,-17.0){$w$}
\drawedge[ELside=r](n1,n11){$b$}
\drawedge[AHnb=0](n1,n0){$a$}
\drawedge[AHnb=0](n00,n11){$a$}

\drawoval(110,-12,40,20,12)
\node(n2)(94.0,-07.0){$v'$}
\node(n22)(94.0,-17.0){$w'$}
\drawedge[AHnb=0](n2,n22){$a$}
\end{picture}
\caption{\small The labeled graphs $\Gamma$ and $\Delta$ in Case (1)}\label{fig: xb differenciation}
\end{figure}
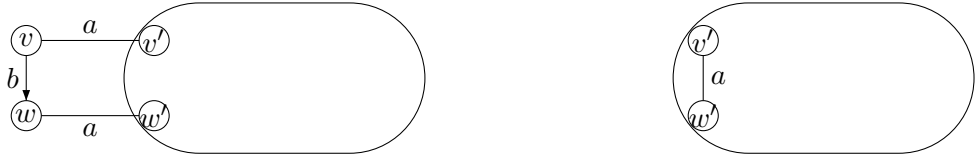

In Case (1), deleting vertices $v$ and $w$ and the edges adjacent to them, and adding an $a$-edge $e'$ between the $v'$ and $w'$ yields a $\PSL$-cyclically reduced graph $\Delta$ with combinatorial type $\bm\tau  + \bm\kappa_{3}$, where $\bm\kappa_{3} =  (-2,-1,-1,0,0)$. If $i$ is the rank of the isolated $a$-edge $e'$ in $\Delta$, this move is called \emph{$(\kappa_3,v, w,i,-)$} if $v' < w'$ and \emph{$(\kappa_3, v, w,i,+)$} if $w' < v'$.

Again, the pair $(\Delta,e')$ and the tuple $(\kappa_3, v, w,i,-)$ (resp. $(\kappa_3, v, w,i,+)$) uniquely specify $(\Gamma,e)$ --- provided that $\Delta$ is  \emph{valid for $(\kappa_3,v,w,i)$}, that is, has no vertex $v$ or $w$ and has at least $i$ isolated $a$-edges. Here as for $\lambda_{2,2}$-moves, the same validity condition allows reversing a $(\kappa_3,v, w,i,-)$-move and a $(\kappa_3,v,w,i,+)$-move, yielding different quasi-labeled graphs $\Gamma$.

Thus $\kappa_3$-moves establish a bijection between the set of structures $(\Gamma,e)$ arising from Case (1) with $\Gamma$ labeled, of size $n$ and combinatorial type $\bm\tau$, and the set of structures $(\Delta,e',v,w,\epsilon)$ formed by a labeled $\PSL$-cyclically reduced graph $\Delta$ of combinatorial type $\bm\tau  + \bm\kappa_{3}$, an isolated $a$-edge $e'$ in $\Delta$, distinct integers $v$ and $w$ in $[n]$ and a boolean $\epsilon$ (to distinguish the cases where $(\Delta,e')$ arises from a $(\kappa_3, v, w,i,-)$- or a $(\kappa_3, v, w,i,+)$-move).

No $\kappa_3$-move is defined in Cases (2), (3) and (4). Observe that, in Case (2), a $\lambda_{2,2}$-move is defined.

\begin{remark}
The description of $\kappa_3$-moves may seem unnaturally complicated, refering as it does to the rank of the isolated $a$-edge between $v'$ and $w'$, rather than to these vertices themselves. This encoding allows us to invert a $\kappa_3$-move in a unique fashion as long as we operate on a graph that is valid for it (a property that is determined by its set of vertex labels and its combinatorial type). This is used in a crucial manner in Section~\ref{sec: generic properties}.
\end{remark}

\paragraph{$\exc$ moves}
For completeness, we also introduce a last, exceptional category of moves, which can be applied only to a quasi-labeled version of $\Delta_3$, turning it into a quasi-labeled version of $\Delta_1$. More precisely, if the $b$-edge in $\Delta_3$ goes from vertex $v$ to vertex $w$, the $\exc(w)$-move returns $\Delta_1$ where the only vertex is labeled $v$. This move can be seen as a degenerate version of a $\lambda_{2,2}$-move. To handle this move like the others, it modifies the combinatorial type by the addition of $\bm{exc} = (-1,0,-1,-1,1)$, the difference between the combinatorial types of $\Delta_1$ and $\Delta_3$.

\subsection{Silhouette of a $\PSL$-cyclically reduced graph}\label{sec: silhouette}

In general, several moves can be applied to a quasi-labeled $\PSL$-cyclically reduced graph $\Gamma$. However, the end result of a maximal sequence of $\lambda_3$-, $\lambda_{2,1}$-, $\lambda_{2,2}$- and $\kappa_{3}$- and $\exc$-moves, in any order, yields the same quasi-labeled graph. Indeed, if two distinct moves are possible, then they modify the $\PSL$-cyclically reduced graph in two disjoint areas which therefore do not interfere with one another: whichever move is performed first, the second one can be performed afterwards.
There are however three exceptions.
\begin{itemize}
\item[(1)] If a $b$-triangle carries exactly two $a$-loops, then two $\lambda_{2,1}$-moves are possible and once one is performed, the other cannot be; but then a $\lambda_{2,2}$-move (or, in a degenerate case, an $\exc$-move) is possible and both options lead to the same quasi-labeled graph after that additional move.
  
\item[(2)] If $\Gamma$ is a $b$-triangle with three $a$-loops, three $\lambda_{2,1}$-moves are possible, each preserving a different $b$-edge and yielding a different quasi-labeled version of $\Delta_3$. We decide to keep the $b$-edge with maximal start vertex, say $v$. The next move can only be an $\exc$-move, leading to the $v$-labeled version of $\Delta_1$.

\item[(3)] If $\Gamma$ is a cycle whose edges spell a word of the form $ab^{\epsilon_1}ab^{\epsilon_2}\cdots ab^{\epsilon_n}$ where $\epsilon_i = \pm1$ for each $i$ and $n\ge 2$, then every $b$-edge is isolated and the only possible moves are $\kappa_3$-moves. A maximal sequence of moves then eliminates all the $b$-edges but one, leading to a quasi-labeled version of $\Delta_2$, with the labels corresponding to the $b$-edge of $\Gamma$ that is not eliminated. We decide to keep the $b$-edge with maximal start vertex.
\end{itemize}

With this discussion in mind and if $\Gamma$ is a labeled $\PSL$-cyclically reduced graph, we define the \emph{quasi-silhouette} $\qcore(\Gamma)$ of $\Gamma$ to be the quasi-labeled graph resulting from the application of a maximal sequence of moves, and the \emph{silhouette} of $\Gamma$ to be $\core(\Gamma) = \relab(\qcore(\Gamma))$.

\begin{example}
Consider the (labeled) Stallings graphs in Figure~\ref{fig: stallings graphs}. The first is equal to its own silhouette. The (quasi-)silhouetting operation applied to the other ones lead to graphs shown in Figure~\ref{fig: silhouettes}.
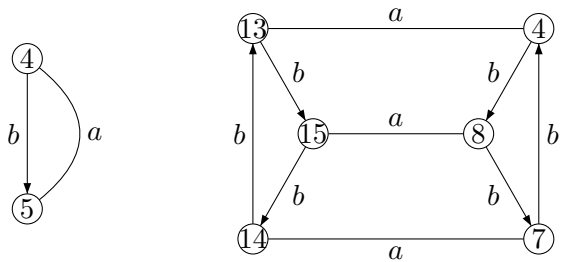
\begin{figure}[htbp]
\centering
\begin{picture}(68,28)(0,-28)
\gasset{Nw=4,Nh=4}
\node(m4)(0,-4){4}
\node(m5)(0,-24.0){5}
\drawedge[ELside=r](m4,m5){$b$}
\drawedge[AHnb=0,curvedepth=7](m4,m5){$a$}
\node(h23)(68.0,-0.0){4}
\node(h12)(68.0,-28.0){7}
\node(h31)(60.0,-14.0){8}
\node(h56)(30.0,0){13}
\node(h46)(30.0,-28.0){14}
\node(h66)(38.0,-14.0){15}
\drawedge[AHnb=0](h56,h23){$a$}
\drawedge[AHnb=0,ELside=r](h31,h66){$a$}
\drawedge[AHnb=0](h12,h46){$a$}
\drawedge[ELside=r](h23,h31){$b$}
\drawedge[ELside=r](h31,h12){$b$}
\drawedge[ELside=r](h12,h23){$b$}
\drawedge(h46,h56){$b$}
\drawedge(h56,h66){$b$}
\drawedge(h66,h46){$b$}
\end{picture}
\caption{\small The (quasi-)silhouettes of the last two (labeled) Stallings graphs in Figure~\ref{fig: stallings graphs}}\label{fig: silhouettes}
\end{figure}
\end{example}

If $\core(\Gamma)$ is $\Delta_1$ or $\Delta_2$, we talk of a \emph{small silhouette}. Otherwise, we talk of a \emph{large silhouette}.
We say that $\Gamma$ is a \emph{silhouette} graph if $\Gamma = \core(\Gamma)$ and $\Gamma \ne \Delta_1, \Delta_2$. The combinatorial type of a silhouette graph is of the form $(n,n/2,0,0,0)$ where $n$ is a multiple of 6. Note that such a graph is the (unrooted) Stallings graph of a free finite index subgroup of $\PSL$, see Section~\ref{sec: rappels}.

The $\PSL$-cyclically reduced subgroups with small silhouette admit the following algebraic characterization.

\begin{proposition}\label{prop: rank and small silhouette}
Let $H$ be a $\PSL$-cyclically reduced subgroup of $\PSL$, with Stallings graph $\Gamma$ and isomorphism type $(\ella, \ellb, r)$. If $\core(\Gamma)$ has isomorphism type $(\ella',\ellb',r')$, then $r = r'$. In particular $\core(\Gamma) = \Delta_1$ (resp. $\Delta_2$) if and only if $r = 0$ (resp. $r = 1$).
\end{proposition}

\proof
Let $\bm\tau$ be the combinatorial type of $\Gamma$. Proposition~\ref{prop: combinatorial vs. isomorphism} shows that the free rank $r$ in the isomorphism type of $H$ is a function of $\bm\tau$; more precisely, if $\bm\tau = (n,\ka,\kb,\ella,\ellb)$, then $6(r-1) = n-2\kb-3\ella-4\ellb = \phi(\bm\tau)$, and we observe that $\phi$ is a linear map.

By construction, $\core(\Gamma)$ is obtained from $\Gamma$ by a succession of $\lambda_3$-, $\lambda_{2,1}$-, $\lambda_{2,2}$-, $\kappa_3$-moves and the $\exc$-move. Each of these moves modifies the combinatorial type by adding to it the vector $\bm\lambda_3$, $\bm\lambda_{2,1}$, $\bm\lambda_{2,2}$, $\bm\kappa_3$ or $\bm{exc}$. Every one of these vectors lies in the kernel of $\phi$, so the free rank component of the isomorphism types of $\Gamma$ and $\core(\Gamma)$ coincide.

It is immediate that this free rank component is 0 for $\Delta_1$, 1 for $\Delta_2$ and $1 + n/6 \ge 2$ for each silhouette graph of size $n$. The proposition follows immediately.
\eop

\subsection{Silhouetting preserves uniformity}\label{sec: randomness preserved}

In Section~\ref{sec: ab-cycles in general}, we will use the fact that taking the silhouette of a uniformly random labeled $\PSL$-reduced graph results in a uniformly random labeled silhouette graph.

More precisely, for $1\leq s\leq n$, let $\cyc(n,s)$ be the set of size $n$ labeled $\PSL$-cyclically reduced graphs whose silhouette has $s$ vertices. Let also 
$\sk(n,s)$ denote the set of quasi-silhouettes of the elements of $\cyc(n,s)$,
$\sk(n,s) = \{\qcore(\Gamma) \mid \Gamma\in\cyc(n,s)\}$. Note that the graphs in $\sk(n,s)$ are labeled by functions into $[n]$ and  that $\sk(s,s)$ is the set of silhouettes of size $s$. We prove the following statement.

\begin{theorem}\label{thm:uniformity}
Let $1\leq s\leq n$. If $\Gamma$ is an element of $\cyc(n,s)$ taken uniformly at random, then $\core(\Gamma)$ is a uniformly random element of $\sk(s,s)$. That is: for any $\Delta,\Delta'\in\sk(s,s)$,
\[
\P(\core(\Gamma)=\Delta)= \P(\core(\Gamma)=\Delta'),
\]
where $\Gamma$ is the $\cyc(n,s)$-valued random variable with uniform distribution.
\end{theorem}

To begin with, we note that there is exactly one possible silhouette of size 1 (resp. size 2), so the statement is trivially true for $s \le 2$. We now assume that $s\ge 3$.

Now we observe that silhouette graphs that differ only by their labeling function have the same probability, as stated in Lemma~\ref{lm:uniformity labels} below. If $\Delta$ is a quasi-labeled graph with labels in $[n]$ and $\sigma$ is a permutation of $[n]$, we let $\sigma(\Delta)$ be the graph obtained from $\Delta$ by relabeling its vertices according to $\sigma$.

\begin{lemma}\label{lm:uniformity labels}
Let $3 \leq s \leq n$, let $\Delta\in\sk(n,s)$ and let $\sigma$ be a permutation of $[n]$, then
\[
\P(\qcore(\Gamma)=\Delta)= \P(\qcore(\Gamma)=\sigma(\Delta)),
\]
where $\Gamma$ is the $\cyc(n,s)$-valued random variable with uniform distribution.
In particular, $\P(\qcore(\Gamma)=\Delta) = \P(\qcore(\Gamma)=\relab(\Delta))$.
\end{lemma}

\proof
Let $\mathfrak{G}$ (resp.  $\mathfrak{G}_\sigma$) be the set of graphs of $\cyc(n,s)$ with quasi-silhouette $\Delta$ (resp. $\sigma(\Delta)$). We have 
\[
\P(\qcore(\Gamma)=\Delta)  = \frac{|\mathfrak{G}|}{|\cyc(n,s)|},
\]
so we only need to establish that $|\mathfrak{G}|=|\mathfrak{G}_\sigma|$.

The moves defined in Section~\ref{sec: moves} do not depend on labels, only on the relative order of labels. So any sequence of moves used to build the quasi-silhouette $\Delta$ from an element $\Gamma\in \mathfrak{G}$, can be applied identically
on $\sigma(\Gamma)$. This shows that $\qcore(\sigma(\Gamma)) = \sigma(\Delta)$ and, in particular, $\sigma(\Gamma)\in\mathfrak{G}_\sigma$.

As
$\sigma$ is a bijection on $\cyc(n,s)$, it follows that $|\mathfrak{G}| \leq |\mathfrak{G}_\sigma|$. The reverse inequality holds by symmetry and this concludes the proof.
\eop

Recall, from Section~\ref{sec: moves}, that if a quasi-labeled graph $\Delta$ is valid for a move $m$, then there exists a unique quasi-labeled graph $\Gamma$ such that $m$ transforms $\Gamma$ into $\Delta$. This extends to a sequence of moves $M = m_t\cdots m_1$: we say that $\Delta$ is \emph{valid for $M$} if there exists a quasi-labeled graph $\Gamma$ such that applying $m_t$, then $m_{t-1}$, etc, yields $\Delta$. In that case, $\Gamma$ is uniquely determined by $\Delta$ and $M$ and we write $\Gamma = M\cdot \Delta$.

We now define, given $\Gamma$, a canonical sequence of moves transforming $\Gamma$ into its quasi-silhouette $\qcore(\Gamma)$. More precisely, the \emph{minimal move} on $\Gamma$ is chosen as follows.

\begin{itemize}
\item If $\Gamma$ has a $b$-loop, the minimal move is the unique move of the form $(\lambda_3,v,w)$, where $v$ is the least vertex carrying a $b$-loop.

\item If $\Gamma$ has an $a$-loop but no $b$-loop and $\Gamma \ne \Delta_3$, the minimal move is the unique move of the form $(\lambda_{2,1},v,w)$, $(\lambda_{2,2}, v\to w,w')$ or $(\lambda_{2,2}, v\leftarrow w,w')$ where $v$ is the least vertex carrying an $a$-loop.

\item If $\Gamma$ has no loop and $\Gamma \ne \Delta_2$, the minimal move is the possible move of the form $(\kappa_3, v, w,i,-)$ or $(\kappa_3, v, w,i,+)$ with $v$ minimal.

\item Finally, if $\Gamma = \Delta_3$, the minimal move is $\exc$.
\end{itemize}
The \emph{minimal sequence} of a graph $\Gamma$ is formed by the minimal move on $\Gamma$, followed by the minimal sequence on the resulting quasi-labeled graph. Observe that no minimal move is defined on $\Gamma$ if $\Gamma = \core(\Gamma)$: in that case, the minimal sequence is the empty sequence $\epsilon$. Thus, the minimal sequence of $\Gamma$ describes a particular way to transform $\Gamma$ into $\qcore(\Gamma)$.

Let $\Delta\in\sk(s,s)$ be a size $s$ labeled silhouette graph. Recall that, by Lemma~\ref{lm:uniformity labels}, the probability that $\qcore(\Gamma)$ is $\Delta$ is the same as the probability of any other quasi-silhouette normalizing to $\Delta$. We establish the following statement.

\begin{lemma}\label{lem:minseq}
Let $\Delta\in\sk(s,s)$. Let $\minseq_n(\Delta)$ be the set of all minimal sequences of the elements of $\cyc(n,s)$ whose quasi-silhouette is $\Delta$. For $\Gamma$ taken uniformly at random in $\cyc(n,s)$, we have
\[
\P(\qcore(\Gamma)=\Delta) = \frac{|\minseq_n(\Delta)|}{|\cyc(n,s)|}.
\]
\end{lemma}

\proof
Every graph $\Gamma\in \cyc(n,s)$ such that $\qcore(\Gamma) = \Delta$ has a unique minimal sequence, which belongs to $\minseq_n(\Delta)$. Conversely, as discussed above, every sequence in $\minseq_n(\Delta)$ uniquely determines the graph in $\cyc(n,s)$ from which it starts. 
\eop

The next lemma is the cornerstone of the proof of Theorem~\ref{thm:uniformity}.
 
\begin{lemma}\label{lm:same minseq}
Let $2 < s \le n$ and let $\Delta$ and $\Delta'$ be labeled silhouette graphs of $\sk(s,s)$. Then
$\minseq_n(\Delta) = \minseq_n(\Delta')$.
\end{lemma} 

\proof
Let $M = m_tm_{t-1}\ldots m_1$ be a minimal sequence in $\minseq_n(\Delta)$.
For $0\leq i\leq t$, let $M_i = m_i\ldots m_1$ be the suffix of $M$ of length $i$, with the convention that $M_0$ is the empty sequence. We verify by induction on $i$ that
\begin{itemize}
\item $\Delta'$ is valid for $M_i$,
\item $M_i\cdot \Delta$ and $M_i\cdot \Delta'$ use the same set of labels,
\item $v$ has a $a$-loop (resp. $b$-loop) in $M_i\cdot \Delta$ iff $v$ has  a $a$-loop (resp. $b$-loop) in $M_i\cdot \Delta'$,
\item $v\rightarrow w$ is an isolated $b$-edge in $M_i\cdot \Delta$ iff $v\rightarrow w$ is an isolated $b$-edge in $M_i\cdot \Delta'$,
\item $M_i\cdot \Delta$ and $M_i\cdot \Delta'$ have the same combinatorial type,
\item if $i\geq 1$, $m_i$ is the minimal move of $M_i\cdot \Delta'$.
\end{itemize}
These properties are true for $i=0$, since $M_0\Delta=\Delta$ and $M_0\Delta'=\Delta'$: both are labeled by $[s]$, have no loop and no isolated $b$-edge, and have combinatorial type $(s,s/2,0,0,0)$.

For the induction step, we observe that inverting a move removes and adds loops with the same labels and on the same vertices on both graphs, and it also removes and adds the same isolated $b$-edges. The only difference that might occur is for $\kappa_3$ moves: the $a$-edge in the move's description is specified by its rank in the order on isolated $a$-edges and not by the labels of its vertices, but the combinatorial type of the quasi-labeled graph evolves the same way. This covers all the items to be proved except for the last one, about the minimality of move $m_i$. This last verification follows directly from the definition of a minimal move.

More precisely, we note that, by definition, $m_i$ is the minimal move for $M_i\cdot \Delta$.
\begin{itemize}
\item If $m_i=(\lambda_3,v,w)$ then $M_i\cdot \Delta$ has a $b$-loop at vertex $v$ and $v$ is minimal such vertex. As the vertices with $b$-loops have the same labels in $M_i\cdot \Delta$ and $M_i\cdot \Delta'$, $m_i$ is also the minimal move of $M_i\cdot \Delta'$.

\item If $m_i$ is of the form $(\lambda_{2,1},v,w)$, $(\lambda_{2,2},v\to w,w')$ or $(\lambda_{2,2},v\leftarrow w,w')$-move, then $M_i\cdot\Delta$ has an $a$-loop, no $b$-loop, an isolated $b$-edge from $v$ to $w$ or the reverse (for a $\lambda_{2,2}$-move) and $v$ is the minimal vertex with an $a$-loop. The same holds by induction for $M_i\cdot\Delta'$ and $m_i$ is minimal for $M_i\cdot\Delta'$.

\item Similarly, if $m_i$ is of the form $(\kappa_3, v, w,i,-)$ or $(\kappa_3, v, w,i,+)$, then both $M_i\Delta$ and $M_i\Delta'$ are loop-free and $m_i$ is minimal for $M_i\Delta'$.

\item Finally, if $m_i = \exc$, then $M_i\cdot\Delta = \Delta_3$, $M_i\cdot\Delta'$ has the same combinatorial type and hence $M_i\cdot\Delta' = \Delta_3$.
\end{itemize}
Thus every sequence of $\minseq_n(\Delta)$ is also in $\minseq_n(\Delta')$. By symmetry, these two sets are equal.
\eop

\proofof{Theorem~\ref{thm:uniformity}}
As indicated before, we may assume that $s > 2$. Let $\Delta$, $\Delta'$ be labeled size $s$ silhouette graphs. 
We have
$$\P(\core(\Gamma) = \Delta) = \sum \P(\qcore(\Gamma) = \Xi)$$
where the sum runs over all $\Xi\in\cyc(n,s)$ such that $\relab(\Xi) = \Delta$. There are $n\choose s$ such graphs $\Xi$, so Lemmas~\ref{lm:uniformity labels} and~\ref{lem:minseq} imply that
\begin{align*}
\P(\core(\Gamma) = \Delta) &= {n\choose s} \P(\qcore(\Gamma) = \Delta) \\
&={n\choose s} \frac{|\minseq_n(\Delta)|}{|\cyc(n,s)|}.
\end{align*}
Similarly,
$$\P(\core(\Gamma) = \Delta') ={n\choose s} \frac{|\minseq_n(\Delta')|}{|\cyc(n,s)|}.$$
We conclude since $\minseq_n(\Delta) = \minseq_n(\Delta')$ by Lemma~\ref{lm:same minseq}.
\eopo

\subsection{Size of the silhouette of a $\PSL$-cyclically reduced graph}\label{sec: size silhouette}

Let $\Gamma$ be a $\PSL$-cyclically reduced graph. The computation of $\core(\Gamma)$ shows how vertices are deleted until there are no more loops or isolated $b$-edges, or until the process has led to $\Delta_1$ or $\Delta_2$. Statistically however, only a small number of vertices are deleted. Proposition~\ref{prop: large silhouette} below quantifies this statement. It follows that most vertices and edges of $\Gamma$ are untouched in the reduction to $\core(\Gamma)$. This observation will play an important role in our discussion of generic properties of subgroups, see Section~\ref{sec: generic properties}, especially Theorem~\ref{thm: small cycles}.

\begin{proposition}\label{prop: large silhouette}
If $\Gamma$ is a $\PSL$-cyclically reduced graph of size $n$, the average number of vertices of $\core(\Gamma)$ is greater than $n - 2 n^\frac23 +o(n^\frac23)$.

Moreover, if $\mu> 0$, there exists $0 < \gamma < 1$ such that the probability that $\core(\Gamma)$ has less than $n - (2+\mu) n^\frac23$ vertices is $\O(\gamma^{n^{\frac13}})$.
\end{proposition}

\proof
Let $\Gamma$ be a $\PSL$-cyclically reduced graph with combinatorial type $(n,\ka,\kb,\ella,\ellb)$ and consider the computation of $\core(\Gamma)$, say by the application of the minimal sequence $M$ of $\Gamma$.

The sequence $M$ starts with $\ellb$ $\lambda_3$-moves, which delete $\ellb$ vertices and add $\ellb$ $a$-loops.
After these moves, $M$ has $\ella+\ellb$ $\lambda_2$-moves. Every $\lambda_{2,1}$-move deletes 1 vertex and 1 $a$-loop, and adds one isolated $b$-edge, so there can be at most $\ella$ such moves. Every $\lambda_{2,2}$-moves deletes 2 vertices, 1 $a$-loop and 1 isolated $b$-edge. Therefore these $\ella+\ellb$ $\lambda_2$-moves delete at most $2(\ella+\ellb)$ vertices. They also create at most $\ella+\ellb$ isolated $b$-edges.

Then $M$ has a sequence of $\kappa_3$-moves, of length the number of  isolated $b$-edges, which is at most $\kb+\ella+\ellb$. Each deletes 2 vertices, so together they delete at most $2(\kb+\ella+\ellb)$ vertices. 

Therefore the total number of vertices deleted in the computation of $\core(\Gamma)$ is at most 
$$\ellb + 2(\ella+\ellb) + 2(\kb+\ella+\ellb) = 2\kb + 4\ella + 5\ellb.$$

The statement on average values is then a direct consequence of this computation and of \cite[Proposition 5.3]{2020:BassinoNicaudWeil}.  Moreover, \cite[Theorem 5.5 and Section 7.2]{2020:BassinoNicaudWeil} establish that, for each $\mu > 0$, there exists $0 < \gamma < 1$ such that
\begin{align*}
\P\left(\kb < (1+\frac\mu6) n^\frac23\right) &= 1-\O\left(\gamma^{n^\frac23}\right),\\
\P\left(\ella < 2 n^\frac12\right) &= 1-\O\left(\gamma^{n^\frac12}\right),\\
\P\left(\ellb < 2 n^\frac13\right) &= 1-\O\left(\gamma^{n^\frac13}\right).
\end{align*}
For $n$ large enough, $2n^\frac13 < 2n^\frac12 < \frac\mu{15}n^{\frac23}$. So
$$\P\left(2\kb + 4\ella + 5\ellb < (2+\mu) n^\frac23\right) = 1-\O\left(\gamma^{n^\frac13}\right),$$
which concludes the proof.
\eop

\subsection{Silhouetting all finitely generated subgroups of $\PSL$}\label{sec: silhouette all subgroups}

Recall that a rooted graph $(\Gamma,v)$ is not $\PSL$-cyclically reduced if and only $v$ is not adjacent to both an $a$- and a $b$-edge, and is the only such vertex. In particular, adding a loop at $v$ labeled by the missing letter(s) yields a $\PSL$-cyclically reduced graph which we denote by $\Gamma^\textsf{o}$. The only case where two loops are added is when $\Gamma$ is the graph with one vertex and no edge (the Stallings graph of the trivial subgroup), which we do not need to consider in the rest of this section.

To deal with labeled $\PSL$-reduced rooted graphs, we introduce a new move, which consists in removing the root and adding the missing loop (if there is a missing loop). More precisely, if $(\Gamma,v)$ is a labeled rooted $\PSL$-reduced graph, let $\alpha = a$ (resp. $b$) if $v$ is not adjacent to an $a$- (resp. $b$-) edge and $\alpha = 0$ if $v$ is adjacent to both an $a$- and a $b$-edge (that is: $\Gamma$ is $\PSL$-cyclically reduced). The \emph{$(\unroot,\alpha,v)$-move} then turns $(\Gamma,v)$ to the labeled $\PSL$-cyclically reduced graph $\Delta = \Gamma^\textsf{o}$ if $\alpha \ne 0$ and $\Delta = \Gamma$ if $\alpha = 0$. Note that exactly one $\unroot$-move is defined on a given $(\Gamma,v)$.

We say that a quasi-labeled $\PSL$-cyclically reduced graph $\Delta$ is \emph{valid} for $(\unroot,\alpha,v)$ if $\Delta$ has a vertex labeled $v$ and if either $\alpha=0$ or $\alpha \ne 0$ and $\Delta$ has an $\alpha$-loop at $v$.

We now define the (quasi-)silhouette of a labeled $\PSL$-reduced rooted graph $(\Gamma,v)$ to be the (quasi-)silhouette of the graph obtained after applying the appropriate $\unroot$-move.

If $1\leq s\leq n$, we let $\rooted(n,s)$ be the set of size $n$ labeled $\PSL$-reduced rooted graphs whose silhouette has $s$ vertices. Observe that 
the set $\sk(n,s)$ of quasi-silhouettes of the elements of $\cyc(n,s)$, is also the set of quasi-silhouettes of the elements of $\rooted(n,s)$. A statement parallel to Theorem~\ref{thm:uniformity} holds for rooted graphs.

\begin{theorem}\label{thm:uniformity rooted}
Let $1\leq s\leq n$. If $(\Gamma,v)$ is an element of $\rooted(n,s)$ taken uniformly at random, then $\core(\Gamma,v)$ is a uniformly random element of $\sk(s,s)$. That is: for any $\Delta,\Delta'\in\sk(s,s)$,
\[
\P(\core((\Gamma,v))=\Delta)= \P(\core((\Gamma,v))=\Delta'),
\]
where $(\Gamma,v)$ is the $\rooted(n,s)$-valued random variable with uniform distribution.
\end{theorem}

\proof
We first observe that Lemma~\ref{lm:uniformity labels} can be readily extended to $\rooted(n,s)$. Then we set the $\unroot$-move as the minimal move on a rooted graph. The \emph{minimal sequence} of a rooted graph consists therefore in an $\unroot$-move, followed by the minimal sequence of the resulting $\PSL$-cyclically reduced graph. There is no difficulty in extending Lemma~\ref{lm:same minseq} to rooted graphs: we consider the sequence $M_{t-1}$ starting at the second move (that is, removing the first move $m_t$, which is an $\unroot$-move); the induction in the proof of the lemma ensures that $M_1$ is the minimal sequence of $M_1\cdot\Delta$ and $M_1\cdot\Delta'$, and that both have the same $a$- and $b$-loops. As the first move $m_t$ is of the $(\unroot,\alpha,v)$-move that is valid for $M_1\cdot\Delta$, it is also valid for $M_1\cdot\Delta'$ (and minimal since an \unroot-move is always minimal). The proof is then completed as in Theorem~\ref{thm:uniformity}.
\eop

We also establish a result similar to Proposition~\ref{prop: large silhouette}.

\begin{proposition}\label{prop: large silhouette rooted}
If $(\Gamma,v)$ is a random labeled rooted $\PSL$-reduced graph of size $n$, the average number of vertices of $\core(\Gamma,v)$ is greater than $n - 2 n^\frac23 +o(n^\frac23)$.

Moreover, if $\mu> 0$, there exists $0 < \gamma < 1$ such that the probability that $\core((\Gamma,v))$ has less than $n - (2+\mu) n^\frac23$ vertices is $\O(\gamma^{n^{\frac13}})$.
\end{proposition}

\proof
Let $\cyc(n)$ and $\rooted(n)$ denote the set of labeled $\PSL$-cyclically reduced graphs and $\PSL$-reduced rooted graphs, respectively.

Let $(\Gamma,v)\in \rooted(n)$ be a $\PSL$-reduced rooted graph and let $\Delta = \unroot(\Gamma,v)$ be the $\PSL$-cyclically reduced graph $\Delta$ obtained after applying the $\unroot$-move. If $(\Gamma,v)$ has combinatorial type $(n,\ka,\kb,\ella,\ellb)$, then $\Delta$ has type $(n,\ka,\kb,\ella',\ellb')$ with $(\ell'_2,\ell'_3) = (\ell_2,\ell_3)$, $(\ell_2+1,\ell_3)$ or $(\ell_2,\ell_3+1)$.

Since we need to manipulate simultaneously random variables defined on different probability spaces (namely the uniform distributions on $\cyc(n)$ and $\rooted(n)$), we introduce the following notation.
\begin{itemize}
\item $\Ella$ (resp. $\Ellb$, $\Kb$) is the random variable defined on $\cyc(n)$ such that $\Ella(\Delta)$ (resp. $\Ellb(\Delta)$, $\Kb(\Delta)$) is the number of $a$-loops (resp. $b$-loops, isolated $b$-edges) in $\Delta$.
\item $\Ella'$ (resp. $\Ellb'$, $\Kb'$) is the random variable defined on $\rooted(n)$ such that $\Ella'(\Gamma,v)$ (resp. $\Ellb'(\Gamma,v)$, $\Kb'(\Gamma,v)$) is the number of $a$-loops (resp. $b$-loops, isolated $b$-edges) in $\Delta = \unroot(\Gamma,v)$ --- that is, $\Ella' = \Ella\circ\unroot$, $\Ellb' = \Ellb\circ\unroot$, $\Kb' = \Kb\circ\unroot$.
\end{itemize}

We have
$$\P(\Ella'(\Gamma,v)=i) = \frac{|\{(\Gamma,v)\in\rooted(n) : \Ella'(\Gamma,v) =i\}|}{|\rooted(n)|}.$$
Observe that a given size $n$ labeled $\PSL$-cyclically reduced graph $\Delta$ is produced by an $\unroot$-move starting from exactly $n + \Ella(\Delta) + \Ellb(\Delta)$ labeled rooted $\PSL$-reduced graphs. If $n > 1$, $n + \Ella(\Delta) + \Ellb(\Delta) \le 2n$, so
\[
|\{(\Gamma,v)\in\rooted(n) : \Ella'(\Gamma,v) =i\}| \enspace\leq\enspace 2n \, |\{\Delta\in\cyc(n) : \Ella(\Delta) =i\}|, 
\]
and since $|\rooted(n)|\geq n |\cyc(n)|$, we have
\[
\P(\Ella'(\Gamma,v)=i) \enspace\leq\enspace \frac{2\, |\{\Delta\in\cyc(n) : \Ella(\Delta) =i\}|}{|\cyc(n)|}
\enspace=\enspace 2 \, \P(\Ella(\Delta) = i).
\]
Similarly we have
\[
\P(\Kb'(\Gamma,v)=i) \enspace\leq\enspace 2 \, \P(\Kb(\Delta)=i)\enspace\text{ and }\enspace
\P(\Ellb'(\Gamma,v)=i) \enspace\leq\enspace  2 \,\P(\Ellb(\Delta)=i).
\]
It follows that we can reproduce the proof of Proposition~\ref{prop: large silhouette}, multiplying the error terms by $2$, which does not change the final statement, expressed using the $\O$ notation. For instance,
\[
\P\left(\Ellb'(\Gamma,v) \geq 2 n^{\frac13}\right) \enspace\leq\enspace 2\,\P\left(\Ellb(\Delta)\geq 2 n^{\frac13}\right) \enspace=\enspace \O\left(\gamma^{n^{1/3}}\right)
\]
for a well chosen value of $\gamma$.
\eop

\section{Counting subgroups by isomorphism and by combinatorial type}\label{sec: counting}

Our aim in this section is to count subgroups of a given size, under some additional constraint: with a fixed isomorphism type or with a fixed combinatorial type. Since subgroups are uniquely represented by their Stallings graph, \emph{i.e.}, by a rooted $\PSL$-reduced graph, this is equivalent to counting these graphs.

It turns out that, for each $n\ge 1$, there are exactly $n!$ distinct labelings of an $n$-vertex rooted $\PSL$-reduced graph, see \textit{e.g.}, \cite[Section 3.1]{2020:BassinoNicaudWeil}. Note that there is no such easy correlation between the number of labeled and unlabeled cyclically reduced graphs, as counting is perturbed by the number of symmetries.

Thus our task reduces to counting labeled $\PSL$-rooted graphs. It further reduces to counting labeled $\PSL$-cyclically reduced graphs, as we explain below.

\subsection{Reduction to the count of labeled $\PSL$-cyclically reduced graphs}\label{sec: reduction to cyclically reduced}

If $\bm\tau = (n,\ka,\kb,\ella,\ellb)$ is a tuple of integers, we let $H(\bm\tau)$ (resp. $L(\bm\tau)$, $s(\bm\tau)$) be the number of subgroups (resp. labeled rooted $\PSL$-reduced graphs, labeled $\PSL$-cyclically reduced graphs) of combinatorial type $\bm\tau$.

\begin{example}\label{ex: S coefficients small n}
In view of Example~\ref{ex: small n}, the non-zero values of $H$, $L$ and $s$ for tuples $(n,\ka,\kb,\ella,\ellb)$ where $n = 1,2$ are as follows:

$\bullet$\enspace $H(\bm\tau)  = L(\bm\tau) = 1$ for $\bm\tau = (1,0,0,1,1), (1,0,0,1,0), (1,0,0,0,1)$ and $s(\bm\tau) = 1$ for $\bm\tau = (1,0,0,1,1)$;

$\bullet$\enspace $L(\bm\tau) = 4$ and $H(\bm\tau) = s(\bm\tau) = 2$ for $\bm\tau = (2,1,1,0,0), (2,0,1,2,0)$; 

$\bullet$\enspace $L(\bm\tau) = 2$ and $H(\bm\tau) = s(\bm\tau) = 1$ for $\bm\tau = (2,1,0,0,2)$; 

$\bullet$\enspace $L(\bm\tau) = 2$ and $H(\bm\tau) = 1$ for $\bm\tau = (2,0,1,1,0), (2,1,0,0,1)$; 
\end{example}

We first establish the connection between the parameters $H(\bm\tau)$, $L(\bm\tau)$ and $s(\bm\tau)$.

\begin{proposition}\label{prop: counting subgroups by cb type}
Let $\bm\tau = (n,\ka,\kb,\ella,\ellb)$ be a combinatorial type with $n\ge 2$. The numbers $H(\bm\tau)$, $L(\bm\tau)$ and $s(\bm\tau)$, respectively of subgroups, labeled rooted $\PSL$-reduced graphs and labeled $\PSL$-cyclically reduced graphs of combinatorial type $\bm\tau$ are related as follows.
$$L(\bm\tau) = n\,s(n,\ka,\kb,\ella,\ellb) + (\ella+1)\,s(n,\ka,\kb,\ella+1,\ellb) + (\ellb+1)\,s(n,\ka,\kb,\ella,\ellb+1) $$
$$H(\bm\tau) = \frac1{n!}\,L(\bm\tau).$$
\end{proposition}

\proof 
Let $(\Gamma,v)$ be a rooted $\PSL$-reduced graph with $n\ge 2$ vertices, such that $\Gamma$ is not $\PSL$-cyclically reduced. Then $v$ is adjacent to an $a$-edge but no $b$-edge, or the opposite. Adding a $b$-loop at $v$ in the first case, an $a$-loop in the second case, yields a rooted $\PSL$-cyclically reduced graph $(\Gamma',v)$. Conversely, if $\Gamma'$ is $\PSL$-cyclically reduced, we get rooted $\PSL$-reduced graphs either by rooting $\Gamma'$ at any one of its vertices, or by rooting $\Gamma'$ at a vertex that carries a loop and deleting that loop. The first equality follows directly.

The second equality follows from the first since a rooted size $n$ $\PSL$-reduced graph has $n!$ distinct labelings.
\eop

Based on Proposition~\ref{prop: combinatorial vs. isomorphism}, which relates the isomorphism type and the combinatorial type of a subgroup, we get the following statement.

\begin{proposition}\label{prop: counting subgroups by isom type}
Let $\bm\sigma = (\ella,\ellb,r)$ be an isomorphism type and let $\ka = \frac12(n-\ella)$.

The number of $\PSL$-cyclically reduced subgroups of size $n$ and isomorphism type $\bm\sigma$ is $n\,s(n,\ka,\kb,\ella,\ellb)$, where $\kb = \frac12(n - 3\ella - 4\ellb - 6r + 6)$.

The number of non-$\PSL$-cyclically reduced subgroups of size $n$ and isomorphism type $\bm\sigma$, where the base vertex is adjacent to an $a$-edge, is $(\ellb + 1)\,s(n,\ka,\kb',\ella,\ellb+1)$, where $\kb' = \frac12(n - 3\ella - 4\ellb - 6r + 2)$.

The number of non-$\PSL$-cyclically reduced subgroups of size $n$ and isomorphism type $\bm\sigma$, where the base vertex is adjacent to a $b$-edge, is $(\ella + 1)\,s(n,\ka,\kb'',\ella+1,\ellb)$, where $\kb'' = \frac12(n - 3\ella - 4\ellb - 6r + 4)$.
\end{proposition}

Propositions~\ref{prop: counting subgroups by cb type} and~\ref{prop: counting subgroups by isom type} effectively reduce the counting of subgroups to the counting of labeled $\PSL$-cyclically reduced graphs of a given combinatorial type, which is investigated in Section~\ref{sec: couting cyclically reduced} below.

\subsection{Counting $\PSL$-cyclically reduced graphs}\label{sec: couting cyclically reduced}

Let $S$ be the multivariate exponential generating series (EGS) of labeled $\PSL$-cyclically reduced graphs, where the different variables account for the components of the combinatorial type, namely
$$S(z,x_2,x_3,y_2,y_3) = \sum_{n,\ka,\kb,\ella,\ellb} \frac{s(n,\ka,\kb,\ella,\ellb)}{n!}\ z^nx_2^{\ka}x_3^{\kb} y_2^{\ella}y_3^{\ellb}.$$
For a series $T$ over these variables, it is convenient to denote by $[z^nx_2^{\ka} x_3^{\kb} y_2^{\ella} y_b^{\ellb}]T$ the coefficient of $z^nx_2^{\ka} x_3^{\kb} y_2^{\ella} y_b^{\ellb}$ in $T$.

We use the so-called \emph{pointing} construction \cite[Theorem II.3]{2009:FlajoletSedgewick} to produce recurrence relations for the $s(\bm\tau)$.  If $u$ is one of $x_2$, $x_3$, $y_2$, $y_3$, let $S_u$ denote the partial derivative $\frac\partial{\partial u}S$. For instance, if $u = x_3$,
$$S_{x_3}(z,x_2,x_3,y_2,y_3) = \sum_{n,\ka,\kb,\ella,\ellb} \kb \left([z^nx_2^{\ka}x_3^{\kb} y_2^{\ella}y_3^{\ellb}]S\right) z^nx_2^{\ka}x_3^{\kb-1} y_2^{\ella}y_3^{\ellb}.$$
In particular, the coefficient $[z^nx_2^{\ka}x_3^{\kb} y_2^{\ella}y_3^{\ellb}](x_3S_{x_3})$
counts the number of structures obtained by pointing (\textit{i.e.}, distinguishing) an isolated $b$-edge
in a labeled $\PSL$-cyclically reduced graph of type $(n,\ka,\kb,\ella,\ellb)$ --- that is: this coefficient is equal to the number of pairs of the form $(\Gamma,e)$ where $\Gamma$ is a labeled $\PSL$-cyclically reduced graph of type $(n,\ka,\kb,\ella,\ellb)$ and $e$ is an isolated $b$-edge in $\Gamma$, divided by $n!$.

An analogous reasoning relates the coefficients of $x_2S_{x_2}$ (resp. $y_2S_{y_2}$, $y_3S_{y_3}$) with the couting of pairs of the form $(\Gamma,e)$ where $e$ is an isolated $a$-edge (resp. a $b$-loop, an $a$-loop).

In the definition of $\lambda_3$-, $\lambda_2$- and $\kappa_3$-moves (Section~\ref{sec: moves}), we described bijections which translate into formulas on the partial derivatives of $S$ as follows.

Let $\bm\tau = (n,\ka,\kb,\ella,\ellb)$ be a combinatorial type with $n\ge 2$ and $\ellb > 0$. The discussion of $\lambda_3$-moves established a bijection between the set of pairs $(\Gamma,\ell)$ of a labeled $\PSL$-cyclically reduced graph $\Gamma$ with combinatorial type $\bm\tau$ with a designated $b$-loop $\ell$, and the set of triples $(\Delta,\ell',v)$ where $\Delta$ is a labeled $\PSL$-cyclically reduced graph of combinatorial type $\bm\tau + \bm\lambda_3 = (n-1,\ka-1,\kb,\ella+1,\ellb-1)$, $\ell'$ is a designated $a$-loop in $\Delta$ and $v$ is an integer in $[n]$. It follows that
\begin{align*}
n!\ [z^nx_2^{\ka}x_3^{\kb} y_2^{\ella}y_3^{\ellb}]\left(y_3S_{y_3}\right) &= n\ \left((n-1)! [z^{n-1}x_2^{\ka-1}x_3^{\kb} y_2^{\ella+1}y_3^{\ellb-1}]\left(y_2S_{y_2}\right)\right) \\
[z^nx_2^{\ka}x_3^{\kb} y_2^{\ella}y_3^{\ellb}]\left(y_3S_{y_3}\right) &= [z^{n-1}x_2^{\ka-1}x_3^{\kb} y_2^{\ella+1}y_3^{\ellb-1}]\left(y_2S_{y_2}\right) \\
&= [z^nx_2^{\ka}x_3^{\kb} y_2^{\ella}y_3^{\ellb}]\left(\frac{zx_2y_2S_{y_2}}{y_2}\right).
\end{align*}
This takes care of the coefficients of $y_3S_{y_3}$ corresponding to Case (1) (that is, where $n > 1$). In Case (2), $\Gamma = \Delta_1$, there is only one $b$-loop to point, so the corresponding EGS is $zy_2y_3$. Therefore we have
\begin{align}
y_3S_{y_3} &= \frac{zx_2y_3}{y_2}\,(y_2S_{y_2}) + zy_2y_3,\nonumber\\
S_{y_3} &= zx_2\,S_{y_2} + zy_2.\label{eq: Stb}
\end{align}

Similarly, the different cases in the discussion of $\lambda_2$-moves yield the following equalities:
\begin{align}
y_2S_{y_2} &= \frac{zy_2}{x_3}\,(x_3S_{x_3}) + 2z^2x_2x_3\,(y_2S_{y_2}) + z^2x_3y_2^2 + zy_2y_3 \nonumber\\
S_{y_2} &= z\,S_{x_3} + 2z^2x_2x_3\,S_{y_2} + z^2x_3y_2 + zy_3.\label{eq: Sta}
\end{align}

Finally, the analysis of $\kappa_3$-moves shows that
\begin{equation}
S_{x_3} = 2z^2x_2^2\,S_{x_2} + 2z^2x_2y_2\,S_{y_2} + z^2x_2 + z^2y_2^2.\label{eq: Sy}
\end{equation}

\subsubsection{Recurrence relations}\label{sec: induction relations}

We now use Equations~(\ref{eq: Stb}), (\ref{eq: Sta}) and (\ref{eq: Sy}) to compute the coefficients $s(n,\ka,\kb,\ella,\ellb)$, with all the arguments non-negative, $n\ge 3$ and $n = 2\ka + \ella \ge 2\kb + \ellb$ (see Example~\ref{ex: small n} for $n < 3$).

Equation~(\ref{eq: Stb}) states that:
\begin{align*}
\sum \frac{\ellb}{n!} s(n,\ka,\kb,\ella,\ellb)\ z^nx_2^{\ka}&x_3^{\kb}y_2^{\ella} y_3^{\ellb-1} \\
&=  zy_2  + \sum \frac{\ella}{n!} s(n,\ka,\kb,\ella,\ellb)\ z^{n+1}x_2^{\ka+1}x_3^{\kb}y_2^{\ella-1}y_3^{\ellb}.
\end{align*}
For $n \ge 3$, $\ellb \ge 1$, considering the coefficient of $z^nx_2^{\ka}x_3^{\kb} y_2^{\ella}y_3^{\ellb-1}$ in this equation, yields the equality
\begin{align}
\frac{\ellb}{n!}\ s(n,\ka,\kb,\ella,\ellb) &= \frac{\ella + 1}{(n-1)!}\ s(n-1,\ka-1,\kb,\ella+1,\ellb-1),\nonumber \\
s(n,\ka,\kb,\ella,\ellb) &= \frac{n(\ella + 1)}{\ellb}\ s(n-1,\ka-1,\kb,\ella+1,\ellb-1).\label{eq: fromStb}
\end{align}

Now turn to Equation~(\ref{eq: Sta}):
\begin{align*}
\sum \frac{\ella}{n!} s(n,\ka,\kb,\ella,\ellb)\ z^nx_2^{\ka}x_3^{\kb}y_2^{\ella-1} & y_3^{\ellb}  =
%
zy_3 + z^2x_3y_2\\
&+ \sum \frac{\kb}{n!} s(n,\ka,\kb,\ella,\ellb)\ z^{n+1}x_2^{\ka}x_3^{\kb-1}y_2^{\ella}y_3^{\ellb}\\
&+ \sum 2\frac{\ella}{n!}\ s(n,\ka,\kb,\ella,\ellb)\ z^{n+2}x_2^{\ka+1}x_3^{\kb+1}y_2^{\ella-1}y_3^{\ellb}.
\end{align*}
For $n \ge 3$, $\ella \ge 1$, considering the coefficient of $z^nx_2^{\ka}x_3^{\kb} y_2^{\ella-1}y_3^{\ellb}$ in this equation, yields 
\begin{align}
\frac{\ella}{n!}\ s(n,\ka,\kb,\ella,\ellb) \enspace=\enspace & \frac{\kb + 1}{(n-1)!}\ s(n-1,\ka,\kb+1,\ella-1,\ellb) \nonumber \\
& + 2\frac{\ella}{(n-2)!}\ s(n-2,\ka-1,\kb-1,\ella,\ellb),\nonumber \\
s(n,\ka,\kb,\ella,\ellb) \enspace=\enspace &\frac{n(\kb + 1)}{\ella}\ s(n-1,\ka,\kb+1,\ella-1,\ellb)\nonumber \\
& + 2n(n-1)\ s(n-2,\ka-1,\kb-1,\ella,\ellb).\label{eq: fromSta}
\end{align}

Finally Equation~(\ref{eq: Sy}) reads as follows:
\begin{align*}
\sum \frac{\kb}{n!} s(n,\ka,\kb,\ella,\ellb)  &z^nx_2^{\ka}x_3^{\kb-1}y_2^{\ella}y_3^{\ellb} \enspace= 
%
z^2x_2 + z^2y_2^2 \nonumber \\
&+ 2\sum \frac{\ka}{n!}\ s(n,\ka,\kb,\ella,\ellb)\ z^{n+2}x_2^{\ka+1}x_3^{\kb}y_2^{\ella}y_3^{\ellb} \\
&+ 2\sum \frac{\ella}{n!}\ s(n,\ka,\kb,\ella,\ellb)\ z^{n+2}x_2^{\ka+1}x_3^{\kb}y_2^{\ella}y_3^{\ellb}.\nonumber
\end{align*}
For $n \ge 3$, $\kb \ge 1$, considering the coefficient of $z^nx_2^{\ka}x_3^{\kb-1} y_2^{\ella}y_3^{\ellb}$ in this equation, yields 
\begin{align}
\frac{\kb}{n!}\ s(n,\ka,\kb,\ella,\ellb) \enspace=\enspace & 2\ \frac{\ka - 1}{(n-2)!}\ s(n-2,\ka-1,\kb-1,\ella,\ellb) \nonumber \\
& + 2\ \frac{\ella}{(n-2)!}\ s(n-2,\ka-1,\kb-1,\ella,\ellb),\nonumber \\
s(n,\ka,\kb,\ella,\ellb) \enspace=\enspace & 2\ \frac{n(n-1)(\ka - 1)}{\kb}\ s(n-2,\ka-1,\kb-1,\ella,\ellb)\nonumber \\
& + 2\ \frac{n(n-1)\ella}{\kb}\ s(n-2,\ka-1,\kb-1,\ella,\ellb).\label{eq: fromSy}
\end{align}

\subsubsection{The base cases and the number of (labeled) silhouette graphs}\label{sec: base cases}

We can use Equations~(\ref{eq: fromStb}), (\ref{eq: fromSta}) and~(\ref{eq: fromSy}) to compute the coefficient $s(n,\ka,\kb,\ella,\ellb)$, where $n \ge 3$: if one of $\kb$, $\ella$ or $\ellb$ is greater than zero, every application of these equations reduces the first argument of the coefficients to compute by 1 or 2. We are therefore reduced to computing coefficients of the form $s(n,\ka,\kb,\ella,\ellb)$ where $n \le 2$, and this was done in Example~\ref{ex: S coefficients small n}, or of the form $s(n,\ka,0,0,0)$, that is, the number of size $n$ silhouette graphs.

The results of  \cite[Proposition 8.18 and Appendix A.4]{2020:BassinoNicaudWeil}  can be used to compute the latter numbers and their asymptotic equivalent (see also the computation by Stothers \cite{1978:Stothers-MathsComp} of the number of finite index, free subgroups of $\PSL$, that is, of subgroups having a silhouette Stallings graph).

\begin{proposition}
Let $T_2$ (resp. $T_3$, $\tilde g$) be given, for $n\ge 1$, by
\begin{align*}
T_2(2n) = \frac{(2n)!}{2^n\, n!} = \prod_{1\le i\le n} (2i-1), \quad
& T_3(3n) = \frac{(3n)!}{3^n\, n!} = \prod_{1\le i\le n} (3i-1)(3i-2) 
\end{align*}
and 
$\tilde g(6n) = T_2(6n)\,T_3(6n).$
Then the number $s(6n,3n,0,0,0)$ of size $6n$ labeled silhouette graphs satisfies $s(0,0,0,0,0) = 0$ and  for $n\ge 1$,
\begin{align*}
s(6n,3n,0,0,0) &= \tilde g(6n) - \sum_{m=1}^{n-1} \tilde g(6m)\,s(6(n-m),3(n-m),0,0,0).
\end{align*}
Moreover $s(6n,3n,0,0,0) \sim \sqrt 6\,\exp\left(7n\log n - 7(1-\log 6)n\right)$.
\end{proposition}

\proof
The first statements in the proposition are taken directly from \cite[Appendix A.4]{2020:BassinoNicaudWeil}. Regarding the last statement, the proof of \cite[Proposition 8.18]{2020:BassinoNicaudWeil} gives an asymptotic equivalent of the EGS $G^\fffi$ of labeled silhouette graphs, namely:
$$[z^{6n}]G^\fffi \sim \frac1{\sqrt 2\pi n} \exp\left(n\log n - (1-\log 6)n\right).$$
Since $s(6n,3n,0,0,0) = (6n)! [z^{6n}]G^\fffi$ and $(6n)! \sim \sqrt{12 \pi n} \exp\left(-6n + 6n(\log(6n))\right)$ --- by Stirling's formula, we have
$$s(6n,3n,0,0,0)  \sim \sqrt 6 \exp\left(7 n\log n - 7(1-\log6)n\right), $$
as announced.
\eop

\section{Random generation of subgroups of $\PSL$}\label{sec: random generation}

As we saw in Example~\ref{ex: S coefficients small n}, there are exactly four size 1 subgroups, with pairwise distinct combinatorial and isomorphism type: the trivial subgroup, the subgroups generated by $a$ and $b$, respectively, and $\PSL$ itself. We now concentrate on generating subgroups of size at least 2, and we assume that the parameters $L(\bm\tau)$ and $s(\bm\tau)$ have been pre-computed for all types of sufficient size.

Like in Section~\ref{sec: counting}, generating uniformly at random a subgroup of a given combinatorial or isomorphism type reduces to randomly generating a labeled rooted $\PSL$-reduced graph of a given combinatorial type and, before that, to randomly generating a labeled $\PSL$-cyclically reduced graph of a given type.

We start with the particular case of labeled silhouette graphs, then proceed to the   general case of labeled $\PSL$-cyclically reduced graphs and, finally, to labeled rooted $\PSL$-reduced graphs.

\subsection{Random labeled silhouette graphs}\label{sec: drawing a silhouette graph}

Let $n$ be a positive multiple of 6. The procedure to generate a size $n$ labeled silhouette graph is simple and well known (see \cite{2020:BassinoNicaudWeil} for instance).

First randomly generate a fixpoint-free permutation on $[n]$ of order 2 and another of order 3. This is done in a standard fashion, by applying a random permutation to the permutation $(1\ 2)(3\ 4)\dots(n-1\ n)$, and another random permutation to $(1\ 2\ 3)(4\ 5\ 6)\dots(n-2\ n-1\ n)$. These two fixpoint-free permutations may determine a disconnected graph, but the proof of  \cite[Proposition 8.18]{2020:BassinoNicaudWeil} shows that this happens with vanishing probability (precisely: $\frac5{36}n\inv + o(n\inv)$). Therefore a rejection algorithm (if the pair of permutations determines a disconnected graph, toss it and draw another pair) produces a silhouette graph after $k$ iterations, with $\E(k) \sim 1$.

\subsection{Random $\PSL$-cyclically reduced graphs}\label{sec: drawing a PSL cyclically reduced graph}

We now consider a tuple $\bm\tau = (n,\ka,\kb,\ella,\ellb)$ such that $s(\bm\tau) > 0$, and we show how we can draw uniformly at random a labeled $\PSL$-cyclically reduced graph with type $\bm\tau$.

We saw in Section~\ref{sec: silhouette} how to construct the (quasi-)silhouette of a $\PSL$-cyclically reduced graph by a sequence of $\lambda_3$-, $\lambda_{2,1}$-, $\lambda_{2,2}$- and $\kappa_3$-moves, possibly followed by an application of the $\exc$-move, each move modifying the combinatorial type by a particular vector. Unwrapping this sequence of moves yields a random generation algorithm which we now explain in more detail.

We construct a finite stochastic path $\path(\bm\tau)$ between types as follows. The edges of this path are labeled by symbols taken in the set $\{\lambda_3, \lambda_{2,1}, \lambda_{2,2}, \kappa_3, \exc\}$.
\begin{itemize}
\item[(1)] If $\bm\tau = (1,0,0,1,1)$, $\bm\tau = (2,1,1,0,0)$ or $\bm\tau$ is of the form $(6s,3s,0,0,0)$ for some $s \ge 1$ (that is, if $\bm\tau$ is the type of $\Delta_1$, $\Delta_2$ or a silhouette graph), we let $\path(\bm\tau)$ be the empty path at $\bm\tau$.
\item[(2)] If $\bm\tau = (2,0,1,2,0)$ (the type of $\Delta_3$), we let $\path(\bm\tau)$ be the length 1 path labeled $\exc$, from $\bm\tau$ to $(1,0,0,1,1) = \bm\tau + \bm{exc}$.
\item[(3)] If $n \ge 2$ and $\ellb > 0$ (the case where we can apply a $\lambda_3$-move), we let $\path(\bm\tau)$ consist in a $\lambda_3$-labeled edge from $\bm\tau$ to $\bm\tau + \bm\lambda_3$, followed by $\path(\bm\tau + \bm\lambda_3)$.

\item[(4)] If $n \ge 2$, no $\lambda_3$-move is possible ($\ellb = 0$) and $\ella > 0$, we have (see Equation~(\ref{eq: fromSta}))
\begin{align*}
s(n,\ka,\kb,\ella,0) = \frac{n(\kb + 1)}{\ella}\ &s(n-1,\ka,\kb+1,\ella-1,0) \\
&+ 2n(n-1)\ s(n-2,\ka-1,\kb-1,\ella,0).
\end{align*}
We choose one of the two summands at random, along the distribution given by the pair $(\frac{n(\kb+1)}{\ella}\ s(n-1,\ka,\kb+1,\ella-1,0), 2n(n-1)\ s(n-2,\ka-1,\kb-1,\ella,0))$. Depending on the outcome of this random choice, we let $\path(\bm\tau)$ consist in either a $\lambda_{2,1}$-labeled edge from $\bm\tau$ to $\bm\tau + \bm\lambda_{2,1}$, followed by $\path(\bm\tau + \bm\lambda_{2,1})$; or a $\lambda_{2,2}$-labeled edge from $\bm\tau$ to $\bm\tau + \bm\lambda_{2,2}$, followed by $\path(\bm\tau + \bm\lambda_{2,2})$.

\item[(4)] If $n > 2$, no $\lambda_3$- or $\lambda_2$-move is possible ($\ella = \ellb = 0$) and a $\kappa_3$-move is possible ($\kb > 0$), we let $\path(\bm\tau)$ consist in an $\kappa_3$-labeled edge from $\bm\tau$ to $\bm\tau + \bm\kappa_3$, followed by $\path(\bm\tau + \bm\kappa_3)$.
\end{itemize}
The label of $\path(\bm\tau)$ is a word $\pathlabel(\bm\tau)$ in the language
$$\lambda_3^*(\lambda_{2,1}+\lambda_{2,2})^*\kappa_3^*\,(1+\exc).$$
By construction, the last type along $\path(\bm\tau)$ is the combinatorial type of $\core(\Gamma)$, where $\Gamma$ is any $\PSL$-cyclically reduced graph of type $\bm\tau$.

Let $\bm\tau_0 = \bm\tau, \bm\tau_1, \ldots, \bm\tau_q$ be the types along $\path(\bm\tau)$ and let $m_1\cdots m_q = \pathlabel(\bm\tau)$. We build a sequence $(\Gamma_i)_{0\le i \le q}$ of labeled $\PSL$-cyclically reduced graphs with combinatorial type $\bm\tau_i$, respectively. The promised randomly chosen labeled $\PSL$-cyclically reduced graph is $\Gamma_0$.

In the large silhouette case (where $\bm\tau_q = (6s,3s,0,0,0)$ for some $s \ge 1$), we draw a labeled silhouette graph $\Gamma_q$ with parameters $\bm\tau_q$ uniformly at random following the procedure in Section~\ref{sec: drawing a silhouette graph}. In the small silhouette case (where $\bm\tau_q$ is the combinatorial type of either $\Delta_1$ or $\Delta_2$), we draw $\Gamma_q$ uniformly at random among the (very few) labeled graphs of type $\bm\tau_q$, see Examples~\ref{ex: small n} and~\ref{ex: S coefficients small n}. Suppose that we have built $\Gamma_i$ for some $i > 0$. 

If $m_i = \lambda_3$, we choose uniformly at random an $a$-loop in $\Gamma_i$ (say, at vertex $w$). Then $\Gamma_{i-1}$ is obtained from $\Gamma_i$ by removing this $a$-loop, adding a new vertex $v$ (with a new label), an $a$-edge between $v$ and $w$, and a $b$-loop at $v$.

If $m_i = \lambda_{2,1}$, we choose uniformly at random an isolated $b$-edge in $\Gamma_i$ (say, from vertex $w$ to vertex $w'$). Then $\Gamma_{i-1}$ is obtained from $\Gamma_i$ by completing this $b$-edge into a $b$-triangle visiting a new vertex $v$ (with a new label), and adding an $a$-loop at $v$.

If $m_i = \lambda_{2,2}$, we choose uniformly at random an $a$-loop in $\Gamma_i$ (say, at vertex $w'$). Then $\Gamma_{i-1}$ is obtained from $\Gamma_i$ by removing that $a$-loop, adding new vertices $v$ and $w$ (with new labels), an isolated $a$-edge between $w$ and $w'$, an isolated $b$-edge between $v$ and $w$, whose orientation (from $v$ to $w$ or from $w$ to $v$) is chosen at random, and an $a$-loop at $v$.

Finally, if $m_j = \kappa_3$, we choose at random an isolated $a$-edge in $\Gamma_i$ (say, between vertices $v'$ and $w'$).  Then $\Gamma_{i-1}$ is obtained from $\Gamma_i$ by removing that isolated $a$-edge, adding new vertices $v$ and $w$ (with new labels) and $a$-edges between $v$ and $v'$, and between $w$ and $w'$, and adding an isolated $b$-edge between $v$ and $w$, whose orientation (from $v$ to $w$ or from $w$ to $v$) is chosen at random.

To summarize: The algorithm to generate a labeled $\PSL$-cyclically reduced graph of type $\bm\tau$ consists in
\begin{itemize}
\item [(1)] computing $\path(\bm\tau)$;
\item[(2)] drawing uniformly at random a labeled $\PSL$-cyclically reduced graph with type the last vertex of $\path(\bm\tau)$; 
\item[(3)] working back along $\path(\bm\tau)$ to generate uniformly at random $\PSL$-cyclically reduced graphs with type the vertices of $\path(\bm\tau)$;
\item[(4)] outputting the graph which corresponds to the first vertex of $\path(\bm\tau)$, that is, a $\PSL$-cyclically reduced graph of type $\bm\tau$.
\end{itemize}

\begin{remark}
The procedure outlined above lazily glosses over the question of labeling because the random generation process described in Section~\ref{sec: drawing a PSL reduced graph}, which eventually calls for forgetting vertex labels, is blind to it. If, however, one wants to randomly generate a \emph{labeled} $\PSL$-cyclically reduced graph, the most expedient is to proceed as above, and then relabel vertices with a random permutation of $[n]$, as is commonly done in random sampling of combinatorial structures (see for instance~\cite[footnote p.12]{1994:FlajoletZimmermannVanCutsem}).
\end{remark}

\subsection{Random subgroups of $\PSL$}\label{sec: drawing a PSL reduced graph}

Let $\bm\tau = (n,\ka,\kb,\ella,\ellb)$ be a combinatorial type. The formula for the number $L(\bm\tau)$ of labeled rooted $\PSL$-reduced graphs of type $\bm\tau$, in Proposition~\ref{prop: counting subgroups by cb type} p. \pageref{prop: counting subgroups by cb type}, suggests the following algorithm to draw uniformly at random a labeled rooted $\PSL$-reduced graph of combinatorial type $\bm\tau$.
\begin{itemize}
\item[(1)] Draw an integer $1 \le p \le L(\bm\tau)$ uniformly at random.
\item[(2)] If $p \le n\,s(\bm\tau)$ and $q$ is the quotient of $p$ by $s(\bm\tau)$ (so that $0\le q < n$), draw uniformly at random a labeled $\PSL$-cyclically reduced graph with combinatorial type $\bm\tau$ and root it at vertex $q+1$.
\item[(3)] If $n\,s(\bm\tau) < p \le n\,s(\bm\tau) + (\ella+1)\,s(n,\ka,\kb,\ella+1,\ellb)$ and $q$ is the quotient of $p - n\,s(\bm\tau)$ by $s(n,\ka,\kb,\ella+1,\ellb)$ (so that $0\le q \le \ella$), draw uniformly at random a labeled $\PSL$-cyclically reduced graph with combinatorial type $(n,\ka,\kb,\ella+1,\ellb)$ (as in Section~\ref{sec: drawing a PSL cyclically reduced graph}), delete the $(q+1)$st $a$-loop (following the order of vertex labels) and root the graph at the vertex where that loop used to be.
\item[(4)] If $n\,s(\bm\tau) + (\ella+1)\,s(n,\ka,\kb,\ella+1,\ellb) < p$ and $q$ is the quotient of $p - n\,s(\bm\tau) - (\ella+1)\,s(n,\ka,\kb,\ella+1,\ellb)$ by $s(n,\ka,\kb,\ella,\ellb+1)$ (so that $0\le q \le \ellb$), draw uniformly at random a labeled $\PSL$-cyclically reduced graph with combinatorial type $(n,\ka,\kb,\ella,\ellb+1)$ (as in Section~\ref{sec: drawing a PSL cyclically reduced graph}), delete the $(q+1)$st $b$-loop (following the order of vertex labels) and root the graph at the vertex where that loop used to be.
\end{itemize}

To draw uniformly at random a subgroup of combinatorial type $\bm\tau$, we first draw a labeled rooted $\PSL$-reduced graph of type $\bm\tau$, and then forget the labeling.

\begin{remark}\label{rk: drawing cyclically reduced by cb type}
To draw uniformly at random a $\PSL$-cyclically reduced subgroup of combinatorial type $\bm\tau$, the algorithm is modified as follows: in step (1), one draws an integer $p$ between 1 and $n\,s(\bm\tau)$; one then applies only step (2).
\end{remark}

Now consider an isomorphism type $\bm\sigma = (\ella,\ellb,r)$. Let $\ka = \frac12(n-\ella)$, $\kb = \frac12(n - 3\ella - 4\ellb - 6r + 6)$, $\kb' = \frac12(n - 3\ella - 4\ellb - 6r + 2)$ and $\kb'' = \frac12(n - 3\ella - 4\ellb - 6r + 4)$.

Proposition~\ref{prop: counting subgroups by isom type} p. \pageref{prop: counting subgroups by isom type} suggests the following algorithm to draw uniformly at random a subgroup of size $n$ and isomorphism type $\bm\sigma$.

\begin{itemize}
\item[(1)] Draw uniformly at random an integer $p$ between 1 and
$$n\,s(n,\ka,\kb,\ella,\ellb) + (\ellb + 1)\,s(n,\ka,\kb',\ella,\ellb+1) + (\ella + 1)\,s(n,\ka,\kb'',\ella+1,\ellb).$$
\item[(2)] If $p \le n\,s(n,\ka,\kb,\ella,\ellb)$, draw uniformly at random a labeled rooted $\PSL$-cyclically reduced graph with combinatorial type $(n,\ka,\kb,\ella,\ellb)$.
\item[(3)] If $n\,s(\bm\tau) < p \le n\,s(\bm\tau) + (\ellb+1)\,s(n,\ka,\kb',\ella,\ellb+1)$ and $q$ is the quotient of $p - n\,s(\bm\tau)$ by $s(n,\ka,\kb',\ella,\ellb+1)$ (so that $0\le q \le \ellb$), draw uniformly at random a labeled $\PSL$-cyclically reduced graph with combinatorial type $(n,\ka,\kb',\ella,\ellb+1)$, delete the $q+1$st $b$-loop (following the order of vertex labels) and root the graph at the vertex where that loop used to be.
\item[(4)] If $n\,s(\bm\tau) + (\ellb+1)\,s(n,\ka,\kb',\ella,\ellb+1) < p$ and $q$ is the quotient of $p - n\,s(\bm\tau) - (\ellb+1)\,s(n,\ka,\kb',\ella,\ellb+1)$ by $s(n,\ka,\kb'',\ella+1,\ellb)$ (so that $0\le q \le \ella$), draw uniformly at random a labeled $\PSL$-cyclically reduced graph with combinatorial type $(n,\ka,\kb,\ella+1,\ellb)$, delete the $q+1$st $a$-loop (following the order of vertex labels) and root the graph at the vertex where that loop used to be.
\end{itemize}

This algorithm can be modified as in Remark~\ref{rk: drawing cyclically reduced by cb type} to draw uniformly at random a $\PSL$-cyclically reduced subgroup of a given isomorphism type.

\begin{remark}
Let us briefly consider the time complexity of this algorithm. The (pre-computation of the $s(\bm\tau)$ for combinatorial types $\bm\tau$ of size at most $n$ takes at most quartic time if the multiplications and additions are performed in time $\O(1)$. As the numbers we consider grow very fast, it is more relevant to consider that encoding a number $N$  takes $\lceil\log_2 N\rceil$ bits, and that the basic operations are not performed in constant time anymore. However, even in this more realistic setting, the pre-computation 
of the $s(\bm\tau)$ is still performed in polynomial time.

For a given $\bm\tau$, of size $n$, computing $\pathlabel(\bm\tau)$ and $\path(\bm\tau)$ takes linear time. If the last type along $\path(\bm\tau)$ has size $s > 2$, the rejection algorithm to produce a random silhouette graph of size $s$ takes polynomial time in average (since the expected number of rejects is $o(1)$). Moreover, each step working back along this path updates the labeled graph constructed so far in a bounded amount of time, and the length of this path is at most $n$.

It follows therefore that the average complexity of the algorithms above, to randomly generate a subgroup of $\PSL$, or a $\PSL$-cyclically reduced subgroup, of type $\bm\tau$ is polynomial in $n$.
\end{remark}


\section{Generic properties of subgroups of $\PSL$}\label{sec: generic properties}

We show that, generically, a finitely generated subgroup of $\PSL$ contains parabolic elements (Corollary~\ref{cor: parabolicity}) and fails to be almost malnormal (Theorem~\ref{thm: negligibility}).

Recall that the \emph{parabolic elements} of $\PSL$ are the conjugates of non-trivial powers of $ab$. A subgroup $H$ of $\PSL$ is said to be \emph{non-parabolic} if it contains no parabolic element.  Also, $H$ is \emph{almost malnormal} if, for every $x \not\in H$, $H \cap H^x$ is finite. It is \emph{malnormal} if each of these intersections is trivial: malnormality coincides with almost malnormality if $H$ is torsion-free (\textit{e.g.} a free subgroup of $\PSL$). 

The proofs of the two results, on parabolic elements and on malnormality, are linked. Indeed, we show that a necessary condition for a subgroup $H \le \PSL$ to be almost malnormal is that $\Gamma(H)$ should not contain a cycle labeled by a power at least 2 of $ab$ (Corollary~\ref{cor: ab-cycles}).

If $H$ has finite index, and in particular if $\Gamma(H)$ is a silhouette graph with at least 2 vertices, both letters $a$ and $b$ label permutations of the vertex set of $\Gamma(H)$ (see Proposition~\ref{prop: charact free and findex}). As a result, $\Gamma(H)$ has a loop at every vertex $v$ labeled $(ab)^{n_v}$ for some $n_v \ge 1$. It follows that $H$ contains parabolic elements. Moreover, if $H$ has size at least 2, then one of the $n_v$ is greater than or equal to 2, so $H$ is not almost malnormal.

Our strategy to deal with general  $\PSL$-reduced rooted (resp.  $\PSL$-cyclically reduced) graphs is the following: we state in Section~\ref{sec: ab-cycles in general} a result on the generic existence of short $ab$-cycles in silhouette graphs (Proposition~\ref{pro:small cycles}). The proof of this result is long and technical, and is deferred to Section~\ref{sec: cycles in silhouette}. We then exploit Proposition~\ref{pro:small cycles} to show that, generically, a $\PSL$-reduced (resp. $\PSL$-cyclically reduced) graph contains an $ab$-cycle of length at least 2 (Theorem~\ref{thm: small cycles}). Our result on parabolic elements follows immediately.

In Section~\ref{sec: negligibility}, we establish the announced characterization of almost malnormality, leading directly to the negligibility of this property.

\subsection{$ab$-cycles in $\PSL$-reduced graphs}\label{sec: ab-cycles in general}

We prove the following theorem.

\begin{theorem}\label{thm: small cycles}
  Let $0 < \alpha < \frac16$. Then a random size $n$ $\PSL$-reduced (resp. $\PSL$-cyclically reduced) graph admits an $ab$-cycle of size at least 2 and and at most $n^\alpha$ with probability $1-\O(n^{-\alpha})$.
\end{theorem}

\begin{corollary}\label{cor: parabolicity}
Let $0 < \alpha < \frac16$. A random size $n$ subgroup (resp. cyclically reduced subgroup) of $\PSL$ is non-parabolic with probability $\O(n^{-\alpha})$.
\end{corollary}

The main ingredients of the proof of Theorem~\ref{thm: small cycles} are:
\begin{itemize}
\item Propositions~\ref{prop: large silhouette} and~\ref{prop: large silhouette rooted}, which state that there exists $0 < \gamma < 1$ such that a size $n$ $\PSL$-cyclically reduced (resp. $\PSL$-reduced rooted) graph has a silhouette with size at least $n-3n^{\frac{2}{3}}$ with probability $1 - \O(\gamma^{n^{\frac{2}{3}}})$ (for the uniform distribution on $\PSL$-cyclically reduced or $\PSL$-reduced graphs).

\item Theorems~\ref{thm:uniformity} and~\ref{thm:uniformity rooted} which state that, under appropriate size constraints, taking the silhouette of a $\PSL$-cyclically reduced (resp. $\PSL$-reduced) graph chosen uniformly at random yields a uniformly random silhouette graph.

\item Proposition~\ref{pro:small cycles} below, which states that, if $0 < \alpha < \frac16$, then, with probability $1 - \O(n^{-\alpha})$, a silhouette graph of size $n$ has an $ab$-cycle of size $m$, with $2\le m \le n^{\alpha}$ (for the uniform distribution on silhouette graphs of a given size).

\end{itemize}

We note that, besides its importance in this proof, Proposition~\ref{pro:small cycles} is of independent interest, as it deals with the probability of the presence of short cycles in certain permutation groups, see the discussion in the introduction.

\begin{proposition}\label{pro:small cycles}
  Let $0 < \alpha < \frac16$ and let $n > 0$ be a multiple of 6. Then a random size $n$ silhouette graph admits an $ab$-cycle of size at least 2 and at most $n^\alpha$ with probability $1-\O(n^{-\alpha})$.
\end{proposition}

\proofof{Theorem~\ref{thm: small cycles}}
Recall that $\cyc(n)$ and $\rooted(n)$ denote the sets of size $n$ labeled $\PSL$-cyclically reduced graphs and $\PSL$-reduced rooted graphs, respectively. We first deal with $\PSL$-cyclically reduced graphs.

To begin with let $0 < \alpha < \frac16$. Let $\probap_\alpha(n)$ be the probability that a graph in $\cyc(n)$ has no $ab$-cycle of size in $[2,n^\alpha]$ and $\probap'_\alpha(n)$ be the probability that an element of $\cyc(n)$ has a silhouette of size at least $n-3n^{\frac{2}{3}}$ and has no $ab$-cycle of length in $[2,n^\alpha]$.

Since by Proposition~\ref{prop: large silhouette}, there exists $0 < \gamma < 1$ such that the probability that an element of $\cyc(n)$ has a silhouette with size less than $n-3n^{\frac{2}{3}}$ is $\O(\gamma^{n^{\frac{2}{3}}})$, we get
$$\probap_\alpha(n) \le \probap'_\alpha(n) + \O(\gamma^{n^{\frac{2}{3}}}).$$

Let now $s \ge n-3n^{\frac23}$ and $\probap'_\alpha(n,s)$ be the probability that an element of $\cyc(n,s)$ (the set of elements of $\cyc(n)$ whose silhouette has size $s$) has no $ab$-cycle of length in $[2,n^\alpha]$. Then
$$\probap'_\alpha(n) = \sum_{s = n-3n^{2/3}}^n \P_n\left(\cyc(n,s)\right)\ \probap'_\alpha(n,s),$$
where $\P_{n}$ is the uniform probability on $\cyc(n)$.

Finally for $\Gamma$ taken uniformly at random in $\cyc(n,s)$ let $\probaq_\alpha(n,s)$ be the probability that both $\core(\Gamma)$ and $\Gamma$ have no $ab$-cycle of size in $[2,n^\alpha]$, and $\probaq'_\alpha(n,s)$ the probability that $\core(\Gamma)$ has such a cycle but $\Gamma$ does not. Then $\probap'_\alpha(n,s) \le \probaq_\alpha(n,s) + \probaq'_\alpha(n,s)$.

Theorem~\ref{thm:uniformity} shows that $\probaq_\alpha(n,s)$ is equal to the probability that a size $s$ silhouette graph has no $ab$-cycle of size in $[2,n^\alpha]$. This is less than or equal to the probability that a size $s$ silhouette graph has no $ab$-cycle of size in $[2,s^\alpha]$ and, according to Proposition~\ref{pro:small cycles}, the latter is $\O(s^{-\alpha})$. Thus, there exists a constant $C$ (independent of $s$ or $n$) such that
$$\probaq_\alpha(n,s) \le Cs^{-\alpha} \le C(n-3n^{\frac23})^{-\alpha} \le C'n^{-\alpha}$$
for another constant $C' > 0$.

About $\probaq'_\alpha(n,s)$,we have
$$\probaq'_\alpha(n,s) = \sum \P_{n,s}\left(\core(\Gamma) = \Delta\right)\ \probaq_\Delta,$$
where $\P_{n,s}$ is the uniform probability on $\cyc(n,s)$, $\probaq_\Delta$ is the probability that a graph in $\cyc(n,s)$ having silhouette $\Delta$ has no $ab$-cycle of size in $[2,n^\alpha]$, and the sum is taken over all size $s$ silhouette graphs $\Delta$ who do have an $ab$-cycle of such a size. 

Let $2 \le \lambda \le n^\alpha$ be the length of an $ab$-cycle in $\Delta$. Let $\Gamma \in \cyc(n,s)$ such that $\core(\Gamma) = \Delta$. In reconstructing $\Gamma$ from $\Delta$ (as in Section~\ref{sec: drawing a PSL cyclically reduced graph}), $a$-edges of $\Delta$ get deleted exactly when undoing a $\kappa_3$-move: an $a$-edge is deleted, two new vertices are added and new $a$- and $b$-edges are added. Since $\Gamma$ has at most $3n^{\frac{2}{3}}$ vertices more that than $\Delta$, there at most $\frac32n^{\frac{2}{3}}$ $\kappa_3$-moves to undo, and the majority of $a$-edges of $\Delta$ are also $a$-edges of $\Gamma$. More precisely, at most $\frac32n^{\frac{2}{3}}$ of the $\frac s2$ $a$-edges of $\Delta$ are not $a$-edges of $\Gamma$.
Thus the probability that an $a$-edge of $\Delta$ fails to be an edge of $\Gamma$ is at most
$$ \frac{3n^{\frac23}}{s}\enspace\le\enspace \frac{3n^{\frac23}}{n-3n^{\frac23}} \enspace=\enspace 3n^{-\frac13}\left(1+\O\left(n^{-\frac13}\right)\right).$$
Then the probability that at least one of the $\lambda$ $a$-edges in the $ab$-cycle under consideration is broken in passing from $\Delta$ to $\Gamma$ is bounded above by $3\lambda n^{-\frac13}\left(1+\O\left(n^{-\frac13}\right)\right)$ and
$$\probaq_\Delta \enspace\le\enspace 3\lambda n^{-\frac13}\left(1+\O\left(n^{-\frac13}\right)\right) \enspace\le\enspace  3n^{\alpha-\frac13}\left(1+\O\left(n^{-\frac13}\right)\right).$$
It follows that
\begin{align*}
\probaq'_\alpha(n,s) &\le 3n^{\alpha-1/3}\left(1+\O\left(n^{-1/3}\right)\right) \\
\probap'_\alpha(n,s) &\le \probaq_\alpha(n,s) + \probaq'_\alpha(n,s) \le C'n^{-\alpha} + 3n^{\alpha-1/3}\left(1+\O\left(n^{-1/3}\right)\right) \le C''n^{-\alpha}
\end{align*}
for some constant $C'' > 0$ (since $\alpha - \frac13 < -\frac16 < -\alpha$). This leads to $\probap'_\alpha(n) \le C''n^{-\alpha}$ and to the announced result $\probap_\alpha(n) = \O(n^{-\alpha})$.

The same proof holds for $\PSL$-reduced rooted graphs, reasoning within $\rooted(n)$ instead of $\cyc(n)$, and using Proposition~\ref{prop: large silhouette rooted} and Theorem~\ref{thm:uniformity rooted} instead of Proposition~\ref{prop: large silhouette} and Theorem~\ref{thm:uniformity}.
\eopo

\subsection{Almost malnormality is negligible}\label{sec: negligibility}

Here we prove the following theorem.

\begin{theorem}\label{thm: negligibility}
Let $0 < \alpha < \frac16$. The probability that a size $n$ subgroup of $\PSL$ is almost malnormal is $\O(n^{-\alpha})$.
\end{theorem}

We start with an elementary characterization of almost malnormality in terms of Stallings graphs given in Proposition~\ref{prop: charact malnormal} below (see \cite[Theorem 7.14]{2007:Markus-Epstein-arXiv} and \cite[Statements 6.7, 6.8, 6.10]{2017:KharlampovichMiasnikovWeil} for more general statements).

An element $g\ne 1$ in $\PSL$ is said to be \emph{cyclically reduced} if it has length 1 or if its normal form starts with $a$ and ends with $b^{\pm1}$, or starts with $b^{\pm1}$ and ends with $a$. It is immediate that every non trivial element of $\PSL$ is conjugated to a cyclically reduced element. Moreover \cite[Theorem IV.2.8]{1977:LyndonSchupp}, two  conjugated cyclically reduced elements are cyclic conjugates of one another (that is: their geodesic representatives are of the form, respectively, $tt'$ and $t't$).

\begin{lemma}\label{lemma: cyclically reduced loops}
Let $H$ be a finitely generated subgroup of $\PSL$ with Stallings graph $(\Gamma(H),v_0)$, let $g \in H$ and let $x$ be of minimal length such that $u = x\inv gx$ is cyclically reduced. Then $\Gamma(H)$ has an $x$-path from $v_0$ to some vertex $v_1$ and a $u$-loop at $v_1$.
\end{lemma}

\proof
If $|u| = 1$, then $xux\inv$ is the normal form of $g$, then $\Gamma(H)$ has a loop at $v_0$ labeled $xux\inv$. It follows that there is a $u$-loop at $v_1$, as announced.

If $|u| > 1$ and without loss of generality, we have $u = au'b^\epsilon$ (for some $\epsilon  = \pm 1$). If $x = 1$, then $\Gamma(H)$ has a $u$-loop at $v_1 = v_0$. If $x \ne 1$, then $x$ cannot end with $a$ or $b^\epsilon$ by minimality of $|x|$, so $x = x'b^{-\epsilon}$. The word $x' b^{-\epsilon}au'b^{-\epsilon} {x'}\inv$ is $\PSL$-reduced and hence is the normal form of $xux\inv = g$. Let $v$ be the vertex reached from $v_0$ reading ${x'}\inv$. Considering the loop at $v_0$ labeled by ${x'}\inv b^{-\epsilon}au'b^{-\epsilon} x'$, we see that $v$ and $v_1$ sit on the same $b$-triangle, and that $u$ labels a loop at $v_1$.
\eop

We can now prove the announced characterization of almost malnormality.

\begin{proposition}\label{prop: charact malnormal}
Let $H$ be a finitely generated subgroup of $\PSL$ and let $(\Gamma,v_0)$ be its Stallings graph. Then
$H$ is almost malnormal if and only if there does not exist distinct vertices $p$ and $q$ in $\Gamma(H)$ and a $\PSL$-reduced word $w$ such that (a) $w$ is not a conjugate of $a$ or $b$ (that is: $w$ has infinite order), and (b) $w$ label loops in $\Gamma(H)$ at $p$ and at $q$.
\end{proposition}

\proof
Let $g,h \in \PSL$, such that $h\not\in H$, $g\in H \cap H^h$ and $g$ has infinite order. Up to conjugating $H$ by an appropriate word, we may assume that $g$ is cyclically reduced (and has length greater than 1). Then $g$ labels a loop in $\Gamma(H)$ at $v_0$.

If $g = y^p$ for some (cyclically) reduced word $y$ and $|p| > 1$, then $g$ labels a loop at $|p|$ distinct vertices along the $g$-loop at $v_0$. We now assume that $g$ is not a proper power.

Since $hgh\inv \in H$, Lemma~\ref{lemma: cyclically reduced loops} shows that $\Gamma(H)$ has an $x$-path from $v_0$ to a vertex $v_1$ and a loop labeled by a cyclically reduced word $u$ such that $hgh\inv = xux\inv$. Moreover, since $hgh\inv$ is conjugated to $g$, the word $u$ is a cyclic conjugate of $g$, say, $g = tt'$ and $u = t't$. Let then $v_2$ be the vertex reached from $v_1$ reading $t'$. There are $tt'$-loops at $v_0$ and $v_2$ and we need to show that $v_0 \ne v_2$.

If $v_0 = v_2$, then $tx\inv \in H$. Moreover
$$hgh\inv = xux\inv = xt'tx\inv = xt\inv tt' tx\inv = (tx\inv)\inv g (tx\inv)$$
so $tx\inv h$ and $g$ commute. By the classical characterization of commuting elements in free products \cite[Theorem 4.5]{1976:MagnusKarrassSolitar}, it follows that either $tx\inv h$ and $g$ sit in the same conjugate of $\langle a\rangle$ or $\langle b\rangle$, or $tx\inv h$ and $g$ sit in the same cyclic subgroup.

The first case is impossible since we assumed that $g$ has infinite order. Therefore there exist $y\in \PSL$, $p,q \in \Z$ such that $g = y^p$ and $tx\inv h = y^q$. We assumed that $g$ is not a proper power, so $p = \pm 1$. Therefore $tx\inv h = g^{pq} \in H$ and hence $h\in H$, a contradiction. %
\eop

In the proof of Theorem~\ref{thm: negligibility}, we use the following consequence of Proposition~\ref{prop: charact malnormal}.

\begin{corollary}\label{cor: ab-cycles}
Let $H$ be a subgroup of $\PSL$. If there exists a $\PSL$-reduced word $w$ which is not a conjugate of $a$ or $b$ and an integer $m\ge 2$ such that $w^m$ labels a loop at some vertex $p$ in the Stallings graph $\Gamma(H)$ but $w$ does not label a loop at $p$, then $H$ fails to be almost malnormal.
\end{corollary}

\proof
The word $w^m$ labels loops at every vertex reached from $p$ reading $w^i$, for $0 \le i < m$ and we conclude by Proposition~\ref{prop: charact malnormal}.
\eop

\proofof{Theorem~\ref{thm: negligibility}}
Theorem~\ref{thm: small cycles} shows that the Stallings graph of a size $n$ random subgroup $H$ contains a simple loop labeled $(ab)^m$ (with $2 \le m \le n^\alpha$) with probability $ 1 - \O(n^{-\alpha})$. The result then follows from Corollary~\ref{cor: ab-cycles}.
\eopo

\subsection{$ab$-cycles in silhouette graphs:  proof of \protect{Proposition~\ref{pro:small cycles}}}\label{sec: cycles in silhouette}

Let $\calG$ be the set of labeled graphs that are disjoint unions of silhouette graphs. Recall that the set of $n$-vertex elements of $\calG$ ($n$ a multiple of 6) is in bijection with the set of pairs $(\sigma_2,\sigma_3)$ of fixpoint-free permutations of $[n]$, the first of order 2 (the $a$-edges) and the second of order 3 (the $b$-edges).

If we fix $\sigma_3$, the set $\calG_{\sigma_3}$ of elements of $\calG$ characterized by pairs of the form $(\sigma_2,\sigma_3)$ has cardinality the number of fixpoint-free, order 2 permutations on $[n]$, namely $(n-1)!!$, where the double factorial of an odd integer $q$ is given by $q!! = q(q-2)(q-4)\dots 1$. Since this value does not depend on $\sigma_3$, it is sufficient to establish the restriction of Proposition~\ref{pro:small cycles} to an arbitrary $\calG_{\sigma_3}$, provided that the constants in the $\O$-notation do not depend on $\sigma_3$.

For convenience, we fix $\sigma_3 = (1\,2\,3)(4\,5\,6)\dots(n-2\,n-1\,n)$, and we write $\calG_n$ instead of $\calG_{\sigma_3}$.
We now concentrate on a particular set of $ab$-cycles. Say that an $ab$-cycle in a graph $G \in \calG$ is \emph{simple} if it visits at most one vertex in each $b$-triangle in $G$. Equivalently, walking along the corresponding $(ab)^m$-labeled path does not require traveling through an $a$-edge in both directions. 

Now let $M = \lfloor n^\alpha\rfloor$ and $\MM = \{2,\dots,M\}$. In the context of this proof, we say that a set is \emph{small} if its cardinality is in $\MM$. Then a \emph{small} $ab$-cycle in an element of $\calG_n$ visits a \emph{small} set of $b$-triangles. We need to prove that, with probability $1-\O(n^{-\alpha})$, a random size $n$ silhouette graph admits a small simple $ab$-cycle.

Let $\calC$ be the set of graphs in $\calG_n$ with at least a small simple $ab$-cycle. Let also $\calS$ be the set of all small sets of $b$-triangles in a graph of $\calG_n$ (recall that all the graphs of $\calG_n$ have the same $b$-triangles). If $I \in \calS$, let $\calC_I$ be the set of graphs in $\calG_n$ containing a simple $ab$-cycle visiting exactly the $b$-triangles in $I$. Then
\[
\calC = \bigcup_{I\in\calS} \calC_I.
\]

Finally, let $\ocalC$ denote the complement of $\calC$ in $\calG_n$ (the elements of $\calG_n$ without a small simple $ab$-cycle). We want to show that $\frac{|\ocalC|}{|\calG_n|} = \frac{|\ocalC|}{(n-1)!!}$ is of the form $\O(n^{-\alpha})$.
By the inclusion-exclusion principle, we have
\begin{align*}
|\ocalC| = |\calG_n| - \sum_{I\in\calS}|\calC_I| + \sum_{\substack{\{I_1,I_2\}\subseteq\calS\\ I_1\ne I_2}} |\calC_{I_1}\cap \calC_{I_2}| - \ldots &= \sum_{\calI\subseteq\calS}(-1)^{|\calI|}\left| \bigcap_{I\in\calI}\calC_I\right|
\\
& = \sum_k (-1)^k \sum_{\substack{\calI\subseteq S\\|\calI| = k}}\left| \bigcap_{I\in\calI}\calC_I\right|.
\end{align*}

Truncating the inclusion-exclusion formula on even or odd cardinalities for $\calI$, yields upper and lower bounds for $|\ocalC|$. For any $\kappa \geq 0$ we have (if $\kappa$ is too large, inequalities become equalities) 
\begin{equation}\label{eq: inclusion-exclusion}
\sum_{k=0}^{2\kappa+1} (-1)^k \sum_{\substack{\calI\subseteq S\\|\calI| = k}}\left| \bigcap_{I\in\calI}\calC_I\right|
\enspace\leq\enspace |\ocalC| \enspace\leq\enspace
 \sum_{k=0}^{2\kappa} (-1)^k \sum_{\substack{\calI\subseteq S\\|\calI| = k}}\left| \bigcap_{I\in\calI}\calC_I\right|
\end{equation}

It turns out to be more convenient to work with tuples of $b$-triangles rather than sets. If $\calI=\{I_1,\ldots,I_k\}$ is a small set of $b$-triangles with cardinality $k$, there are $k!$ tuples
$\bm J = (J_1,\ldots,J_k)\in\calS^k$ such that $\{J_1,\ldots,J_k\}=\calI$, so
\begin{equation}
\left| \bigcap_{I\in\calI}\calC_I\right| = \frac1{k!}\sum_{\substack{(J_1,\ldots,J_k)\in\calS^k\\\{J_1,\ldots,J_k\}=\calI}} \left| \calC_{J_1}\cap\ldots\cap\calC_{J_k}\right|.
\label{eq: given set of triangles}
\end{equation}

Say that two simple $ab$-cycles \emph{overlap} if they visit a same $b$-triangle (of course, on different vertices). We now distinguish the tuples $\bm J = (J_1,\ldots,J_k)$ in Equation~\eqref{eq: given set of triangles} according to the cardinality of their components and to their overlaps. More precisely, if $\bm d = (d_1,\ldots,d_k) \in \MM^k$, we let
\begin{align*}
\nooverlap(\bm d) &= \left\{ \bm J = (J_1,\ldots,J_k)\in\calS^k \mid \forall i,\ |J_i|=d_i\text{ and the $J_i$ are pairwise disjoint}\right\} \\
\overlap(\bm d) &= \Big\{ \bm J = (J_1,\ldots,J_k)\in\calS^k \mid \forall i,\ |J_i|=d_i\text{ and  the $J_i$ are pairwise distinct} \\
&\hskip 8cm\text{but not pairwise disjoint}\Big\}.
\end{align*}

Returning to the summands in the estimation of the inclusion-exclusion bounds (Equation~\eqref{eq: inclusion-exclusion}), we now have, for every $k$,
\begin{equation}\label{eq: Ak Bk}
\sum_{\substack{\calI\subseteq S\\|\calI| = k}}\left| \bigcap_{I\in\calI}\calC_I\right|
\enspace=\enspace \frac1{k!}\left(\underbrace{\sum_{\bm d\in\MM^k}\sum_{\bm J\in\nooverlap(\bm d)}\left| \bigcap_{i=1}^k\calC_{J_i}\right|}_{A_k}+\underbrace{\sum_{\bm d\in\MM^k}\sum_{\bm J\in\overlap(\bm d)}\left| \bigcap_{i=1}^k\calC_{J_i}\right|}_{B_k}\right)
\end{equation}
We now study successively the quantities $A_k$ and $B_k$ in Equation~\eqref{eq: Ak Bk}. If $q\ge 1$, we let $H_q$ denote the partial sum $\sum_{i=1}^q\frac1i$ of the harmonic series. It is well known that $H_q = \log q + \gamma + o(1)$, where $\gamma$ is Euler's constant.

\begin{lemma}\label{lemma:Ak}
Let $n$ be a positive multiple of 6. Let $0 < \alpha < \frac12$, $M = \lfloor n^\alpha\rfloor$, $0 < \beta < \frac12-\alpha$ and $1 \le k \le n^\beta$. Finally, let $\delta=\alpha+\beta$ and let $A_k$ be as in Equation~\eqref{eq: Ak Bk}. Then 
\begin{equation*}
\frac1{(n-1)!!}A_k = (H_M-1)^k\left(1 + \O\left(n^{2\delta-1}\right)\right),
\end{equation*}
uniformly in $k$ (that is: the constants intervening in the $\O$ notation do not depend on $n$ or $k$).
\end{lemma}

\proof
Let $\bm d = ( d_1,\dots, d_k) \in\MM^k$, $ d= d_1+\ldots+ d_k$ and $\bm J = (J_1,\ldots,J_k)\in\nooverlap(\bm d)$. To construct a graph in $\bigcap_i\calC_{J_i}$, that is, a graph in $\calG_n$ with (non-overlapping) simple $ab$-cycles over the sets of $b$-triangles $J_1, \dots, J_k$ respectively, we must
\begin{itemize}
\item select for each $b$-triangle in $J_1, \dots, J_k$ a vertex belonging to the collection of simple $ab$-cycles; there are 3 possibilities for each $b$-triangle, and therefore a total of $3^d$ choices;
\item cyclically order the triangles in each $J_i$; there are $( d_1-1)!\dots( d_k-1)!$ possibilities to do so; note that this second step fully determines the $a$-edges adjacent to the $ d$ vertices chosen in the first step;
\item choose the missing $a$-edges arbitrarily: they connect the $n-2 d$ vertices not yet adjacent to an $a$-edge, and there are $(n-2 d-1)!!$ ways to do so.
\end{itemize}

Thus, for every $\bm J = (J_1,\ldots,J_k)\in\nooverlap(\bm d)$, there are $3^ d( d_1-1)!\cdots( d_k-1)!(n-2 d-1)!!$ graphs in $\bigcap_{i=1}^k\calC_{J_i}$ (which is independent of the choice of $\bm J$ in $\nooverlap(\bm d)$).

Moreover, since a graph in $\calG_n$ has $\frac n3$ $b$-triangles and the components of a tuple in $\nooverlap(\bm d)$ are pairwise disjoint, we have $|\nooverlap(\bm d)| = \binom{n/3}{ d_1,\ldots, d_k,n- d}$. Thus
$$\sum_{\bm J\in\nooverlap(\bm d)}\left| \bigcap_{i=1}^k\calC_{J_i}\right|
= \binom{n/3}{ d_1,\ldots, d_k,n/3- d}3^d( d_1-1)!\cdots( d_k-1)!(n-2 d-1)!!.$$
Note that $\frac{(n-2  d-1)!!}{(n-1)!!}= \prod_{i=0}^{ d-1}\frac{1}{n-1-2i}$ and $\frac{(n/3)! 3^{ d}}{(n/3- d)!}=\prod_{i=0}^{ d-1} (n-3i)$. Therefore
$$\frac1{(n-1)!!}\ \sum_{\bm J\in\nooverlap(\bm d)}\left| \bigcap_{i=1}^k\calC_{J_i}\right| \enspace=\enspace  \frac{Q_d}{ d_1\cdots d_k},\quad\text{ with } Q_d = \prod_{i=0}^{d-1}\frac{n-3i}{n-1-2i}.$$
Now observe that
\begin{align}
Q_d &\enspace\leq\enspace \frac{n}{n-1}\quad\text{and} \label{eq: upper bound of Qell}\\
Q_d &\enspace=\enspace  \prod_{i=0}^{d-1}\left(1-\frac{i-1}{n-1-2i} \right) \enspace\ge\enspace \left(1-\frac{d}{n-2d}\right)_d \enspace=\enspace \exp\left(d\log\left(1-\frac{d}{n-2d}\right)\right) \nonumber \\
&\enspace\ge\enspace 1 + d\log\left(1-\frac{d}{n-2d}\right).\label{eq: lower bound of Qell}
\end{align}
The lower bound in Equation \eqref{eq: lower bound of Qell} for $Q_d$ is a decreasing function of $d$ and the possible values of $d$ satisfy  $1\le d \le kM \le n^{\alpha+\beta} = n^\delta$, so
\begin{equation}\label{eq: U and L}
1 + n^\delta\log\left(1-\frac{n^\delta}{n-2n^\delta}\right) \enspace\le\enspace Q_d \enspace\le\enspace \frac n{n-1}.
\end{equation}
Let $U$ and $L$ be the upper and lower bounds of $Q_d$ in Equation~\eqref{eq: U and L}. Note that $\sum_{d\in\MM^k}\frac1{d_1\cdots d_k} = (H_M-1)^k$. Therefore
$$L (H_M-1)^k \enspace=\enspace L \sum_{\bm d\in\MM^k}\frac1{d_1\dots d_k} \enspace\le\enspace \frac1{(n-1)!!}A_k \enspace\le\enspace  U \sum_{\bm d\in\MM^k}\frac1{d_1\dots d_k} \enspace=\enspace U(H_M-1)^k,$$
which concludes the proof of Lemma~\ref{lemma:Ak} since both $U$ and $L$ are of the form $1 + \O\left(n^{2\delta-1}\right)$.
\eop

\begin{lemma}\label{lemma:Bk}
 Let $n$ be a positive multiple of 6. Let $0 < \alpha < \frac16$, $M = \lfloor n^\alpha\rfloor$, $0 < \beta < \frac16-\alpha$ and $1 \le k \le n^\beta$. Finally, let $\delta = \alpha+\beta$ and let $B_k$ be as in Equation~\eqref{eq: Ak Bk}. Then
\begin{equation*}
\frac1{(n-1)!!} B_k \leq \O\left(n^{6\delta-1}\right) (H_M-1)^k
\end{equation*}
uniformly in $k$.
\end{lemma}

\proof
Let $\bm d = (d_1,\dots,d_k) \in\MM^k$ and $d=d_1+\ldots+d_k$. Since every vertex of a silhouette graph occurs in exactly one $ab$-cycle, a $b$-triangle can occur in at most 3 $ab$-cycles. Let $\bm t =(t_1,t_2,t_3)$ be a tuple of non-negative integers such that $d = t_1+2t_2+3t_3$. We denote by $\overlap_k(\bm d;\bm t)$ the set of elements $\bm J = (J_1,\ldots,J_k)\in\overlap_k(\bm d)$ such that there are $t_1$ (resp. $t_2$, $t_3$) $b$-triangles belonging to exactly 1 (resp. 2, 3) of the components of $\bm J$. We talk of $b$-triangles \emph{of type} 1 (resp. 2, 3).

\begin{lemma}\label{lemma: even}
If $\bm J \in \overlap_k(\bm d;\bm t)$ and $\bigcap_{i=1}^k\calC_{J_i} \ne \emptyset$, then $t_2+t_3$ is even.
\end{lemma}

\proof
Let $G \in \bigcap_{i=1}^k\calC_{J_i}$. Consider the $a$-edges occurring in the $k$ small simple $ab$-cycles visiting, respectively, $J_1,\dots, J_k$. It is convenient at this point to think of the $a$-edges in our $ab$-cycles as matched pairs of half-edges, which we denote $(T,e)$, where $e$ is an $a$-edge and $T$ is a $b$-triangle adjacent to $e$. Since no $a$-edge in a simple $ab$-cycle may connect vertices from the same $b$-triangle, each such $a$-edge $e$ corresponds to a pair of distinct, matched half-edges $(T,e)$ and $(T',e)$.

Even though $a$-edges are undirected (or can be traversed in both directions), considering one in a simple $ab$-cycle uniquely defines a direction and we can talk of matched outgoing and incoming half-edges along a simple $ab$-cycle. Clearly, an $a$-edge is used in only one direction in a given $ab$-cycle, but it can be used in different directions by distinct $ab$-cycles. So we say that a half-edge $(T,e)$ such that $e$ occurs in the union of the $k$ $ab$-cycles under consideration is \emph{outgoing} (resp. \emph{incoming}, \emph{2-way}) if it only occurs as outgoing (resp. it only occurs as incoming, it occurs both as outgoing and incoming).

If $T$ is a $b$-triangle in $\bigcup J_i$, occurring in just one of the $k$ small simple $ab$-cycles, then $T$ is a component of exactly 1 incoming and 1 outgoing half-edges; if $T$ occurs in two of these $ab$-cycles, it is a component of 1 incoming, 1 outgoing and 1 2-way half-edges; finally, if $T$ occurs in three $ab$-cycles, it is a component of 3 2-way half-edges.

The result follows since a 2-way half-edge must be matched with another, distinct, 2-way half-edge, and we have $t_2+3t_3$ such half-edges.
\eop

Let $\bm J = (J_1,\ldots,J_k) \in \overlap_k(\bm d;\bm t)$. Then $|\bigcup_i J_i| = t_1+t_2+t_3$ whereas $d = \sum_i|J_i| = t_1+2t_2+3t_3$.  The overlaps between the $J_i$ determine which $b$-triangles in $\bigcup_i J_i$ are of type 1, 2 or 3. To construct an $n$-vertex graph in $\bigcap_i\calC_{J_i}$, that is, a graph in $\calG_n$ with (overlapping) simple $ab$-cycles over the sets of $b$-triangles $J_1, \dots, J_k$ respectively, we must
\begin{itemize}
\item select for each $b$-triangle of type 1 a vertex belonging to the collection of simple $ab$-cycles; there are $3^{t_1}$ choices. Similarly, for each $b$-triangle of type 2, select two vertices belonging to the collection of simple $ab$-cycles; there are $3^{t_2}$ choices. Note that every vertex of a type 3 $b$-triangle belongs to the collection of simple $ab$-cycles. As discussed in the proof of Lemma~\ref{lemma: even}, the choice of these $t_1+2t_2+3t_3$ vertices in the $ab$-cycles implies the presence of $2t_1+3t_2+3t_3$ vertices along the corresponding $(ab)^{d_i}$-labeled loops.
\item cyclically order the triangles in each $J_i$; there are $(d_1-1)!\dots(d_k-1)!$ possibilities to do so; note that this second step fully determines the $a$-edges adjacent to the $2t_1 + 3t_2 + 3t_3$ vertices determined in the first step;
\item choose the missing $a$-edges arbitrarily: they connect the $n-(2t_1 + 3t_2 + 3t_3)$ vertices not yet selected, and there are $(n-(2t_1 + 3t_2 + 3t_3)-1)!!$ ways to do so.
\end{itemize}

Thus, for every $\bm J = (J_1,\ldots,J_k) \in \overlap_k(\bm d;\bm t)$, we have
$$\left|\bigcap_{i=1}^k\calC_{J_i}\right| \enspace=\enspace 3^{t_1+t_2}(d_1-1)!\cdots(d_k-1)!(n-2t_1 - 3t_2 - 3t_3 -1)!!.$$
Now, to construct a tuple in $\overlap_k(\bm d;\bm t)$, we must 
\begin{itemize}
\item choose $t_i$ $b$-triangles of type $i$ ($i = 1, 2, 3$); there are
$\binom{n/3}{t_1,t_2,t_3,n/3-t_1-t_2-t_3}$ choices;
\item allocate the $t_1+t_2+t_3$ selected $b$-triangles to sets $J_1,\dots,J_k$ respecting multiplicities and the required cardinality of the $J_i$; there are at most  $\binom{t_1+2t_2+3t_3}{d_1,\ldots,d_k} = \binom{d}{d_1,\ldots,d_k}$ choices (this is an upper bound as some choices may be unrealizable or produce cycles that are not simple).
\end{itemize}

Therefore we have
$$\left|\overlap_k(\bm d;\bm t)\right| \enspace\le\enspace \binom{n/3}{t_1,t_2,t_3,n/3-t_1-t_2-t_3}\ \binom{d}{d_1,\ldots,d_k}$$
and $\sum_{\bm J\in\overlap(\bm d;\bm t)}\left| \bigcap_{i=1}^k\calC_{J_i}\right|$ is bounded above by
$$\binom{n/3}{t_1,t_2,t_3,n/3-t_1-t_2-t_3}\ \binom{d}{d_1,\ldots,d_k}3^{t_1+t_2}(d_1-1)!\cdots(d_k-1)!(n-2t_1 - 3t_2 - 3t_3 -1)!!.$$
Dividing by $(n-1)!!$, we get
$$\frac{\sum_{\bm J\in\overlap(\bm d;\bm t)}\left| \bigcap_{i=1}^k\calC_{J_i}\right|}{(n-1)!!} \enspace\le\enspace \frac{d !}{d_1\cdots d_k\ t_1!t_2!t_3!\  3^{t_3}}\ 
\frac{n(n-3)\cdots(n-3(t_1+t_2+t_3-1))}{(n-1)(n-3)\cdots(n-2t_1-3t_2-3t_3+1)}.$$
Let $t=t_1+t_2+t_3$, so that $t_1=d - 2t_2-3t_3$. Recall that $Q_t =  \prod_{i=0}^{t-1}\frac{n-3i}{n-1-2i}\leq \frac{n}{n-1}$ (Equation~\eqref{eq: upper bound of Qell}). Then
\begin{align}
\frac{\sum_{\bm J\in\overlap(\bm d;\bm t)}\left| \bigcap_{i=1}^k\calC_{J_i}\right|}{(n-1)!!}
& \enspace\leq\enspace \frac1{d_1\cdots d_k}\cdot \frac{d!}{(d-2t_2-3t_3)!}
\cdot\frac{Q_t}{\prod_{i=t}^{t_1+\frac32(t_2+t_3)-1}(n-1-2i)} \nonumber\\
& \enspace\leq\enspace \frac{n\ d^{2t_2+3t_3}}{(n-1)\ d_1\cdots d_k\ (n-2t_1-3t_2-3t_3+1)^{\frac12(t_2+t_3)}} \label{eq: step in Bk}
\end{align}
Since $d = t_1 + 2t_2 + 3t_3$ and $d \leq kM \leq n^{\alpha+\beta} = n^\delta$, we have $2t_1+3t_2+3t_3 \leq 2d \leq 2n^\delta$, yielding
\[
\frac{d^{2t_2+3t_3}}{(n-2t_1-3t_2-3t_3+1)^{\frac12(t_2+t_3)}} \leq 
\frac{n^{\delta(2t_2+3t_3)}}{(n-2n^\delta)^{\frac12(t_2+t_3)}}.
\]
Moreover, $\frac12(t_2+t_3)\leq \frac14 d \leq \frac14 n^\delta$ and we have
\begin{align*}
(n-2n^\delta)^{\frac12(t_2+t_3)}
& = n^{\frac12(t_2+t_3)}\exp\left(\left(\frac12(t_2+t_3)\right)\log\left(1-2n^{\delta-1}\right)\right)\\
& \geq n^{\frac12(t_2+t_3)}\exp\left(\frac14 n^\delta \log\left(1-2n^{\delta-1}\right)\right)
\end{align*}
The term under the exponential is
$$\frac14 n^\delta \log\left(1-2n^{\delta-1}\right) \sim -\frac12 n^{3\delta-1},$$
which tends to 0 (since $\delta < \frac16 < \frac13$). Therefore, for $n$ sufficiently large, 
\[
(n-2n^\delta)^{\frac12(t_2+t_3)} \geq \frac12 n^{\frac12(t_2+t_3)}.
\]
Going back to Equation~\eqref{eq: step in Bk}, we have
\begin{equation}\label{eq: 2nd step in Bk}
\frac{\sum_{\bm J\in\overlap(\bm d;\bm t)}\left| \bigcap_{i=1}^k\calC_{J_i}\right|}{(n-1)!!} \leq  \frac{1}{d_1\cdots d_k} \frac{n^{\delta(2t_2+3t_3)}}{n^{\frac12(t_2+t_3)}} = \frac{1}{d_1\cdots d_k} \ n^{\delta(2t_2+3t_3)-\frac12(t_2+t_3)}\\
\end{equation}
Equation~\eqref{eq: 2nd step in Bk} must be summed over $\bm t$ and $\bm d$. Let us first fix $\bm d$ and recall that $\bm t=(t_1,t_2,t_3)$ with $(t_2,t_3)\neq (0,0)$ and $t_2+t_3$ even (Lemma~\ref{lemma: even}).

Consider first the $\bm t$ of the form $(t_1,t_2,0) = (t_1,2s,0)$ ($s\ge 1$). The corresponding subsum of powers of $n$ is bounded above (for $n\ge 2$) by
$$\sum_{s\ge 1}\left(n^{2\delta-\frac12}\right)^{2s} = \frac{n^{4\delta -1}}{1 - n^{4\delta -1}} \le 2 n^{4\delta-1}.$$
Similarly, the subsum of powers of $n$ corresponding to the $\bm t$ of the form $(t_1,0,t_3)$ is bounded above by $2 n^{6\delta-1}$.

Next, the subsum corresponding to the $\bm t$ where both $t_2$ and $t_3$ are non-zero is bounded above by
$$\sum_{t_2\ge 1}\sum_{t_3\ge 1} \left(n^{2\delta-\frac12}\right)^{t_2}\left(n^{3\delta-\frac12}\right)^{t_3} \le \left(\sum_{t_2\ge 1}\left(n^{2\delta-\frac12}\right)^{t_2}\right)\ \left(\sum_{t_3\ge 1}\left(n^{3\delta-\frac12}\right)^{t_3}\right) \le 2 n^{5\delta-1}$$
It follows that, for a fixed tuple $\bm d$,
\[
\frac{\sum_{\bm J\in\overlap(\bm d;\bm t)}\left| \bigcap_{i=1}^k\calC_{J_i}\right|}{(n-1)!!} 
\leq \frac{2(n^{4\delta-1} + n^{5\delta-1} + n^{6\delta-1})}{d_1\cdots d_k}.
\]
Now, summing over $\bm d$, we get
\begin{align*}
\frac{\sum_{\bm J\in\overlap(\bm d;\bm t)}\left| \bigcap_{i=1}^k\calC_{J_i}\right|}{(n-1)!!} 
&\leq 2(n^{4\delta-1} + n^{5\delta-1} + n^{6\delta-1})\sum_{\bm d}\frac1{d_1\cdots d_k} \\
&\le 2(n^{4\delta-1} + n^{5\delta-1} + n^{6\delta-1})(H_M-1)^k
\end{align*}
Since $n^{4\delta-1} + n^{5\delta-1} + n^{6\delta-1} = \O\left(n^{6\delta-1}\right)$, we have, for $n$ large enough,
\[
\frac{\sum_{\bm J\in\overlap(\bm d;\bm t)}\left| \bigcap_{i=1}^k\calC_{J_i}\right|}{(n-1)!!} = (H_M-1)^k \O(n^{6\delta-1}),
\]
thus concluding the proof of Lemma~\ref{lemma:Bk}.
\eop

We can now conclude the proof of Proposition~\ref{pro:small cycles}. 
By Equations~\eqref{eq: inclusion-exclusion} and~\eqref{eq: Ak Bk}, we want to show that, for $\kappa=\lfloor n^\beta\rfloor$ and $\kappa=\lfloor n^\beta\rfloor -1$,
\[
\frac1{(n-1)!!}\ \sum_{k=0}^{\kappa}(-1)^k\sum_{\substack{\calI\subseteq S\\|\calI| = k}}\left| \bigcap_{I\in\calI}\calC_I\right|
= \sum_{k=0}^{\kappa}\frac{(-1)^k}{k!}\frac{A_k}{(n-1)!!} + \sum_{k=0}^{\kappa}\frac{(-1)^k}{k!}\frac{B_k}{(n-1)!!}
\]
is $\O(n^{-\alpha})$. 
By Lemma~\ref{lemma:Ak}, the absolute value of the first sum is bounded above by
$$(1+\O(n^{2\delta-1})) \exp(-(H_M-1)) = e^{1-\gamma}n^{-\alpha}(1+\O(n^{2\delta-1})).$$
And by Lemma~\ref{lemma:Bk}, the absolute value of the second sum is bounded above by
$$\exp(-(H_M-1))\ \O(n^{6\delta-1}) = e^{1-\gamma}n^{-\alpha} \O(n^{6\delta-1}).$$
Thus the whole sum is $\O(n^{-\alpha})$, establishing the expected bound for disjoint unions of silhouette graphs.

Such a union, of size $n$ (a multiple of 6), is connected (and hence silhouette) with probability $1 - \frac56n\inv + o(n\inv)$ by \cite[Proof of Proposition 8.18]{2020:BassinoNicaudWeil}. Thus the probability that a silhouette graph has no small simple $ab$-cycle is, again, $\O(n^{-\alpha})$.

{\small\bibliographystyle{abbrv}
\bibliography{pwbiblio}}


\end{document}